# THE STABLE MANIFOLD THEOREM FOR STOCHASTIC DIFFERENTIAL EQUATIONS*

SALAH-ELDIN A. MOHAMMED
MICHAEL K.R. SCHEUTZOW

March 16, 1998

ABSTRACT. We formulate and prove a *Local Stable Manifold Theorem* for stochastic differential equations (sde's) that are driven by spatial Kunita-type semimartingales with stationary ergodic increments. Both Stratonovich and Itô-type equations are treated. Starting with the existence of a stochastic flow for a sde, we introduce the notion of a hyperbolic stationary trajectory. We prove the existence of invariant random stable and unstable manifolds in the neighborhood of the hyperbolic stationary solution. For Stratonovich sde's, the stable and unstable manifolds are dynamically characterized using forward and backward solutions of the anticipating sde. The proof of the stable manifold theorem is based on Ruelle-Oseledec multiplicative ergodic theory.

## 1. Introduction.

Consider the following Stratonovich and Itô stochastic differential equations (sde's) on $\mathbf{R}^d$:

$$\left. \begin{aligned} d\phi(t) &= \overset{\circ}{F}(\circ dt, \phi(t)), \quad t > s \\ \phi(s) &= x \end{aligned} \right\} \quad (S)$$

$$\left. \begin{aligned} d\phi(t) &= F(dt, \phi(t)), \quad t > s \\ \phi(s) &= x \end{aligned} \right\} \quad (I)$$

defined on a filtered probability space $(\Omega, \mathcal{F}, (\mathcal{F}_s^t)_{s \leq t}, P)$. Equation (S) is driven by a continuous forward-backward spatial semimartingale $\overset{\circ}{F} : \mathbf{R} \times \mathbf{R}^d \times \Omega \to \mathbf{R}^d$, and equation (I) is driven by a continuous forward spatial semimartingale

1991 *Mathematics Subject Classification.* Primary 60H10, 60H20; secondary 60H25, 60H05..

Mohammed's research is supported in part by NSF Grants DMS-9503702, DMS-9703852 and by MSRI, Berkeley, California. Scheutzow's research is supported in part by MSRI, Berkeley, California. Research at MSRI is supported in part by NSF grant DMS-9701755.

$F : \mathbf{R} \times \mathbf{R}^d \times \Omega \to \mathbf{R}^d$ (Kunita [Ku]). Both $\overset{\circ}{F}$ and $F$ have stationary ergodic increments.

It is known that, under suitable regularity conditions on the driving spatial semimartingale $\overset{\circ}{F}$, the sde (S) admits a continuous (forward) stochastic flow $\phi_{s,t} : \mathbf{R}^d \times \Omega \to \mathbf{R}^d$, $-\infty < s \leq t < \infty$ ([Ku]). The inverse flow is denoted by $\phi_{t,s} := \phi_{s,t}^{-1} : \mathbf{R}^d \times \Omega \to \mathbf{R}^d$, $-\infty < s \leq t < \infty$. This flow is generated by Kunita's backward Stratonovich sde

$$\left. \begin{aligned} d\phi(s) &= -\overset{\circ}{F}(\circ \hat{d}s, \phi(s)), \quad s < t \\ \phi(t) &= x. \end{aligned} \right\} \quad (S^-)$$

Similarly, the inverse flow $\phi_{t,s} := \phi_{s,t}^{-1} : \mathbf{R}^d \times \Omega \to \mathbf{R}^d$, $-\infty < s \leq t < \infty$ of the Itô equation (I) solves a backward Itô sde with a suitable correction term ([Ku], p. 117).

In this article, we prove a local stable-manifold theorem for the sde's (S) and (I) under the condition that the driving semimartingales $\overset{\circ}{F}$ and $F$ have stationary ergodic increments. Our main result is Theorems 3.1. It gives a random flow-invariant local splitting of $\mathbf{R}^d$ into stable and unstable differentiable submanifolds in the neighborhood of each hyperbolic stationary solution. The method we use to establish these results is based on a non-linear discrete-time multiplicative ergodic theorem due to Ruelle ([Ru.1], cf. [Ru.2]). Although the article is largely self-contained, familiarity with the arguments in [Ru.1] will sometimes be needed. Key ingredients of this approach are Ruelle-Oseledec integrability conditions which we prove in Lemma 3.1. The proof of this lemma is in turn based on spatial estimates on the flow and its derivatives ([Ku], [M-S.2]). These estimates are stated in Theorem 2.1 for easy reference.

During the past few years, several authors have contributed to the development of the stable-manifold theorem for non-linear sde's . The first successful attempt was carried out by Carverhill for sde's on compact manifolds ([C]). In [C], a stable manifold theorem is obtained in the globally asymptotically stable case where the Lyapunov exponents of the linearized flow are all negative. The general hyperbolic case with positive Lyapunov exponents is not treated in [C]. The work by Boxler ([Bo]) focuses on the existence of a (global) center manifold under small (white) noise. Wanner ([Wa]) deals essentially with the real-noise case. In this paper, we are able to handle a much more general type of noise with a simple and transparent approach. In addition, we expect the method of construction of the unstable manifold to work if the state space is replaced by a finite-dimensional Riemannian manifold, cf. [C].

The multiplicative ergodic theory of linear finite-dimensional systems was initiated by Oseledec in his fundamental work ([O]). An infinite-dimensional stable-manifold theorem for linear stochastic delay equations was developed by Mohammed [Mo.1] in the white-noise case, and by Mohammed and Scheutzow for general semi-martingales with stationary ergodic increments ([M-S.1]).

## 2. Basic setting and preliminary results.

Let $(\Omega, \mathcal{F}, P)$ be a probability space. Let $\theta : \mathbf{R} \times \Omega \to \Omega$ be a $P$-preserving flow on $\Omega$, viz.

(i) $\theta$ is jointly measurable,
(ii) $\theta(t+s, \cdot) = \theta(t, \cdot) \circ \theta(s, \cdot), \quad s, t \in \mathbf{R}$,
(iii) $\theta(0, \cdot) = I_\Omega$, the identity map on $\Omega$,
(iv) $P \circ \theta(t, \cdot)^{-1} = P, \quad t \in \mathbf{R}$.

Denote by $\bar{\mathcal{F}}$ the $P$-completion of $\mathcal{F}$. Let $(\mathcal{F}_s^t : \infty < s \leq t < \infty)$ be a family of sub-$\sigma$-algebras of $\bar{\mathcal{F}}$ satisfying the following conditions:

(i) $\theta(-r, \cdot)(\mathcal{F}_s^t) = \mathcal{F}_{s+r}^{t+r}$ for all $r \in \mathbf{R}, -\infty < s \leq t < \infty$.
(ii) For each $s \in \mathbf{R}$, both $(\Omega, \bar{\mathcal{F}}, (\mathcal{F}_s^{s+u})_{u \geq 0}, P)$ and $(\Omega, \bar{\mathcal{F}}, (\mathcal{F}_{s-u}^s)_{u \geq 0}, P)$ are filtered probability spaces satisfying the usual conditions ([Pr.2]).

A random field $F : \mathbf{R} \times \mathbf{R}^d \times \Omega \to \mathbf{R}^d$ is called a *(continuous forward) spatial semimartingale helix* if it satisfies the following:

(i) For every $s \in \mathbf{R}$, there exists a sure event $\Omega_s \in \mathcal{F}$ such that

$$F(t+s, x, \omega) = F(t, x, \theta(s, \omega)) + F(s, x, \omega)$$

for all $t \in \mathbf{R}$, all $\omega \in \Omega_s$ and all $x \in \mathbf{R}^d$.

(ii) For almost all $\omega \in \Omega$, the mapping $\mathbf{R} \times \mathbf{R}^d \ni (t, x) \mapsto F(t, x, \omega) \in \mathbf{R}^d$ is continuous.

(iii) For any fixed $s \in \mathbf{R}$ and $x \in \mathbf{R}^d$, the process $F(s+t, x, \omega) - F(s, x, \omega), t \geq 0$, is an $(\mathcal{F}_s^{s+t})_{t \geq 0}$-semimartingale.

Similarly, a random field $F : \mathbf{R} \times \mathbf{R}^d \times \Omega \to \mathbf{R}^d$ is called a *continuous backward spatial semimartingale helix* if it satisfies (i) and (ii) and has the property that for fixed $s \in \mathbf{R}$ and $x \in \mathbf{R}^d$, the process $F(s-t, x, \omega) - F(s, x, \omega), t \geq 0$, is an $(\mathcal{F}_{s-t}^s)_{t \geq 0}$-semimartingale.

In (iii) above, it is enough to require that the semimartingale property holds for some fixed $s$ (e.g. $s = 0$); then it will hold automatically for every $s \in \mathbf{R}$ ([A-S], Theorem 14).

Note that a semimartingale helix $F$ always satisfies $F(0, x, \omega) = 0$ for a.a. $\omega \in \Omega$ and all $x \in \mathbf{R}^d$. It is also possible to select a suitable *perfect* version of $F$ such that the helix property (i) holds for *every* $\omega \in \Omega$. See [A-S] for further details, and [Pr.1] for other general properties of semimartingale helices.

Suppose that the continuous forward semimartingale helix $F : \mathbf{R} \times \mathbf{R}^d \times \Omega \to \mathbf{R}^d$ is decomposed as

$$F(t, x) = V(t, x) + M(t, x), \quad t \in \mathbf{R}^+, \quad x \in \mathbf{R}^d,$$

where $V(\cdot, x) := (V^1(\cdot, x), \cdots, V^d(\cdot, x))$ is a continuous bounded variation process, $M(\cdot, x) := (M^1(\cdot, x), \cdots, M^d(\cdot, x))$ is a continuous local martingale with respect to $(\mathcal{F}_0^t)_{t \geq 0}$, and $V(0, x) = M(0, x) = 0$ for each $x \in \mathbf{R}^d$. Let $\langle M^i(\cdot, x), M^j(\cdot, y) \rangle$ be the joint quadratic variation of $M^i(\cdot, x), M^j(\cdot, y)$, for $x, y \in \mathbf{R}^d, 1 \leq i, j \leq d$.

Throughout this paper, assume that $F|[0,\infty)$ has forward local characteristics $(a(t,x,y), b(t,x))$ that satisfy the relations

$$\langle M^i(\cdot, x), M^j(\cdot, y)\rangle(t) = \int_0^t a^{i,j}(u,x,y)\,du, \quad V^i(t,x) = \int_0^t b^i(u,x)\,du$$

for all $1 \leq i,j \leq d$, $0 \leq t \leq T$, and where $a(t,x,y) := (a^{i,j}(t,x,y))_{i,j=1,\cdots,d}$, $b(t,x) := (b^1(t,x), \cdots, b^d(t,x))$. Further measurability properties of the local characteristics are given in ([Ku], pp. 79-85). Note that the local characteristics are uniquely determined by $F$ up to null sets.

In what follows, let $\triangle$ denote the diagonal $\triangle := \{(x,x) : x \in \mathbf{R}^d\}$ in $\mathbf{R}^d \times \mathbf{R}^d$, and let $\triangle^c$ be its complement. The space $\mathbf{R}^d$ carries the usual Euclidean norm $|\cdot|$.

We shall use the notation

$$\alpha := (\alpha_1, \alpha_2, \cdots, \alpha_d), \quad D_x^\alpha := \frac{\partial^{|\alpha|}}{(\partial x_1)^{\alpha_1}\cdots(\partial x_d)^{\alpha_d}}, \quad |\alpha| := \sum_{i=1}^d \alpha_i,$$

for $\alpha_i$ non-negative integers, $i = 1, \cdots, d$.

Following [Ku], we shall say that the spatial forward semimartingale $F$ has forward local characteristics of class $(B_{ub}^{m,\delta}, B_{ub}^{k,\delta})$ for non-negative integers $m, k$ and $\delta \in (0,1]$, if for all $T > 0$, its characteristics satisfy

$$\operatorname*{ess\,sup}_{\omega \in \Omega} \sup_{0 \leq t \leq T} [\|a(t)\tilde{\|}_{m+\delta} + \|b(t)\|_{k+\delta}] < \infty,$$

where

$$\|a(t)\tilde{\|}_{m+\delta} := \sup_{x,y \in \mathbf{R}^d} \frac{|a(t,x,y)|}{(1+|x|)(1+|y|)} + \sum_{1 \leq |\alpha| \leq m} \sup_{x,y \in \mathbf{R}^d} |D_x^\alpha D_y^\alpha a(t,x,y)|$$
$$+ \sum_{|\alpha|=m} \|D_x^\alpha D_y^\alpha a(t,\cdot,\cdot)\hat{\|}_\delta,$$

$$\|b(t)\|_{k+\delta} := \sup_{x \in \mathbf{R}^d} \frac{|b(t,x)|}{(1+|x|)} + \sum_{1 \leq |\alpha| \leq k} \sup_{x \in \mathbf{R}^d} |D_x^\alpha b(t,x)|$$
$$+ \sum_{|\alpha|=k} \sup_{(x,y) \in \triangle^c} \frac{|D_x^\alpha b(t,x) - D_y^\alpha b(t,y)|}{|x-y|^\delta},$$

and

$$\|f\hat{\|}_\delta := \sup\left\{\frac{|f(x,y) - f(x',y) - f(x,y') + f(x',y')|}{|x-x'|^\delta|y-y'|^\delta} \; : \; (x,x'), (y,y') \in \triangle^c\right\}$$

for any $\delta$-Hölder continuous function $f : \mathbf{R}^d \times \mathbf{R}^d \to \mathbf{R}^d$. Similar definitions hold for the backward local characteristics of a backward spatial semimartingale.

Now consider the Stratonovich and Itô stochastic differential equations

$$\left.\begin{array}{l} d\phi(t) = \overset{\circ}{F}(\circ dt, \phi(t)), \quad t > s, \\ \phi(s) = x, \end{array}\right\} \tag{S}$$

$$\left.\begin{array}{l} d\phi(t) = F(dt, \phi(t)), \quad t > s, \\ \phi(s) = x. \end{array}\right\} \tag{I}$$

The sde (S) is driven by a continuous forward(-backward) spatial helix semimartingale $\overset{\circ}{F}(t,x) := (\overset{\circ}{F}_1(t,x), \cdots, \overset{\circ}{F}_d(t,x))$, $x \in \mathbf{R}^d$. In the sde (I), $F$ denotes a spatial continuous *forward* helix semimartingale. It is known that, under suitable regularity hypotheses on the local characteristics of $\overset{\circ}{F}$ (or $F$), the sde's (S) and (I) generate the *same* stochastic flow. Throughout this article, these flows will denoted by the same symbol $\{\phi_{s,t} : s, t \in \mathbf{R}\}$. More precisely, we will need the following hypotheses:

**Hypothesis (ST$(k,\delta)$).**

$\overset{\circ}{F}$ is a continuous spatial helix forward semimartingale with forward local characteristics of class $(B_{ub}^{k+1,\delta}, B_{ub}^{k,\delta})$. The function

$$[0,\infty) \times \mathbf{R}^d \ni (t,x) \mapsto \sum_{j=1}^d \left.\frac{\partial a^{\cdot,j}(t,x,y)}{\partial x_j}\right|_{y=x} \in \mathbf{R}^d$$

belongs to $B_{ub}^{k,\delta}$.

**Hypothesis (ST$^-(k,\delta)$).**

$\overset{\circ}{F}$ is a continuous helix backward semimartingale with backward local characteristics of class $(B_{ub}^{k+1,\delta}, B_{ub}^{k,\delta})$.

**Hypothesis (IT$(k,\delta)$).**

$F : \mathbf{R} \times \mathbf{R}^d \times \Omega \to \mathbf{R}^d$ is a continuous spatial helix forward semimartingale with forward local characteristics of class $(B_{ub}^{k,\delta}, B_{ub}^{k,\delta})$.

The following proposition establishes a relationship between the sde's (S) and (I).

**Proposition 2.1.**

*Suppose the helix semimartingale $\overset{\circ}{F}$ satisfies Hypothesis $(ST(k,\delta))$ for some positive integer $k$ and $\delta \in (0,1]$. Let the following relation hold:*

$$F(t,x,\omega) := \overset{\circ}{F}(t,x,\omega) + \frac{1}{2}\int_0^t \sum_{j=1}^d \left.\frac{\partial a^{\cdot,j}(u,x,y)}{\partial x_j}\right|_{y=x} du, \quad t \in \mathbf{R}, x \in \mathbf{R}^d.$$

*Then $F$ is a helix semimartingale which satisfies Hypothesis $(IT(k,\delta))$. In this case, the sde's (S) and (I) generate the same stochastic flow $\phi_{s,t}$, $s,t \in \mathbf{R}$, on $\mathbf{R}^d$.*

*Proof.*

The assertion of the proposition follows from Theorem 3.4.7 in [Ku], except for the helix property. The helix property of $F$ follows from that of $\overset{\circ}{F}$ and the fact that the $\mathbf{R}^{d \times d}$-valued process $<\overset{\circ}{F}(\cdot,x), \overset{\circ}{F}(\cdot,y)>(t) = <F(\cdot,x), F(\cdot,y)>(t)$, $t \in \mathbf{R}$, is a helix for any $x, y \in \mathbf{R}^d$ ([Pr.1]). □

Proposition 2.1 shows that for a given $k, \delta$, Hypothesis $(ST(k,\delta))$ is stronger than $(IT(k,\delta))$. Although our results will cover both the Stratonovich and Itô cases, the reader may note that the Stratonovich sde (S) allows for a complete and more aesthetically pleasing dynamic characterization of the stochastic flow $\phi_{s,t}$ and its inverse. Indeed, under $(ST(k,\delta))$ and $(ST^-(k,\delta))$, $\phi_{s,t}^{-1}$ solves the backward Stratonovich sde based on $\overset{\circ}{F}$ and hence provides a natural dynamical representation of the local unstable manifold in terms of trajectories of the backward Stratonovich sde. Such a dynamical characterization is not available for the Itô sde (I). See Section 3.

From now on, we will implicitly assume that the spatial semimartingales $\overset{\circ}{F}$ and $F$ are related by the formula in Proposition 2.1. In this context, all our results will be derived under both sets of hypotheses $(ST(k,\delta))$ and $(IT(k,\delta))$, although the conclusions pertain invariably to the generated flow $\phi_{s,t}$.

The following proposition is elementary. Its proof is an easy induction argument using the chain rule.

**Proposition 2.2.**

*Let $f := (f_1, f_2, \cdots, f_d) : \mathbf{R}^d \to \mathbf{R}^d$ be a $C^k$ diffeomorphism for some integer $k \geq 1$. Then, for each $\alpha$ with $1 \leq |\alpha| \leq k$,*

$$D_x^\alpha (f)_i^{-1}(x) = \frac{p_{\alpha,i}(f^{-1}(x))}{|\det[Df(f^{-1}(x))]|^{n_\alpha}}, \quad i = 1, \cdots, d, \tag{1}$$

*for all $x \in \mathbf{R}^d$, and some integer $n_\alpha \geq 1$. In the above identity, $p_{\alpha,i}(y)$ is a polynomial in the partial derivatives of $f$ of order up to $|\alpha|$ evaluated at $y \in \mathbf{R}^d$.*

*Proof.*

We use induction on $\alpha$. For $|\alpha| = 1$, the chain rule gives $Df^{-1}(x) = [Df(f^{-1}(x))]^{-1}$. By Cramer's rule, this implies (1) with $n_\alpha = 1$.

Assume by induction that for some integer $1 \leq n < k$, (1) holds for all $\alpha$ such that $|\alpha| \leq n$ and all $i = 1, \cdots, d$. Take $\alpha$ such that $|\alpha| = n$, and fix $i, j \in \{1, 2, \cdots, d\}$. Taking partial derivatives with respect to $x_j$ in both sides of (1) shows that the right-hand-side of the resulting equation is again of the same form with $\alpha$ replaced by $\tilde{\alpha} := (\tilde{\alpha}_1, \tilde{\alpha}_2, \cdots, \tilde{\alpha}_d)$, where $\tilde{\alpha}_i := \alpha_i + \delta_{i,j}$. This completes the proof of the proposition. □

The next proposition allows the selection of sure $\theta(t, \cdot)$-invariant events in $\mathcal{F}$ from corresponding ones in $\bar{\mathcal{F}}$.

**Proposition 2.3.**
Let $\Omega_1 \in \bar{\mathcal{F}}$ be a sure event such that $\theta(t, \cdot)(\Omega_1) \subseteq \Omega_1$ for all $t \geq 0$. Then there is a sure event $\Omega_2 \in \mathcal{F}$ such that $\Omega_2 \subseteq \Omega_1$ and $\theta(t, \cdot)(\Omega_2) = \Omega_2$ for all $t \in \mathbf{R}$.

*Proof.*

Define $\hat{\Omega}_1 := \bigcap_{k=0}^{\infty} \theta(k, \cdot)(\Omega_1)$. Then $\hat{\Omega}_1$ is a sure event, $\hat{\Omega}_1 \subseteq \Omega_1$ and $\theta(t, \cdot)(\hat{\Omega}_1) = \hat{\Omega}_1$ for all $t \in \mathbf{R}$. Since $\bar{\mathcal{F}}$ is the completion of $\mathcal{F}$, we may pick a sure event $\Omega_0 \subseteq \hat{\Omega}_1$ such that $\Omega_0 \in \mathcal{F}$. Define

$$\Omega_2 := \{\omega : \omega \in \Omega, \theta(t, \omega) \in \Omega_0 \text{ for Lebesgue-a.e } t \in \mathbf{R}\}.$$

Using Fubini's theorem and the $P$-preserving property of $\theta$, it is easy to check that $\Omega_2$ satisfies all the conclusions of the proposition. □

**Theorem 2.1.**
Let $\overset{\circ}{F}$ satisfy Hypothesis $(ST(k, \delta))$ (resp. $F$ satisfies $(IT(k, \delta))$) for some $k \geq 1$ and $\delta \in (0, 1]$. Then there exists a jointly measurable modification of the trajectory random field of (S) (resp. (I)) also denoted by $\{\phi_{s,t}(x) : -\infty < s, t < \infty, x \in \mathbf{R}^d\}$, with the following properties:
If $\phi : \mathbf{R} \times \mathbf{R}^d \times \Omega \to \mathbf{R}^d$ is defined by

$$\phi(t, x, \omega) := \phi_{0,t}(x, \omega), \quad x \in \mathbf{R}^d, \omega \in \Omega, t \in \mathbf{R},$$

then the following is true for all $\omega \in \Omega$:
(i) For each $x \in \mathbf{R}^d$, and $s, t \in \mathbf{R}$, $\phi_{s,t}(x, \omega) = \phi(t - s, x, \theta(s, \omega))$.
(ii) $(\phi, \theta)$ is a perfect cocycle:

$$\phi(t + s, \cdot, \omega) = \phi(t, \cdot, \theta(s, \omega)) \circ \phi(s, \cdot, \omega),$$

for all $s, t \in \mathbf{R}$.
(iii) For each $t \in \mathbf{R}$, $\phi(t, \cdot, \omega) : \mathbf{R}^d \to \mathbf{R}^d$ is a $C^k$ diffeomorphism.

(iv) *The mapping* $\mathbf{R}^2 \ni (s,t) \mapsto \phi_{s,t}(\cdot,\omega) \in \mathit{Diff}^k(\mathbf{R}^d)$ *is continuous, where* $\mathit{Diff}^k(\mathbf{R}^d)$ *denotes the group of all $C^k$ diffeomorphisms of $\mathbf{R}^d$, given the $C^k$-topology.*

(v) *For every* $\epsilon \in (0,\delta), \gamma, \rho, K, T > 0$, *and* $1 \leq |\alpha| \leq k$, *the quantities*

$$\sup_{\substack{0 \leq s,t \leq T, \\ x \in \mathbf{R}^d}} \frac{|\phi_{s,t}(x,\omega)|}{[1+|x|(\log^+|x|)^\gamma]}, \qquad \sup_{\substack{0 \leq s,t \leq T, \\ x \in \mathbf{R}^d}} \frac{|D_x^\alpha \phi_{s,t}(x,\omega)|}{(1+|x|^\gamma)},$$

$$\sup_{x \in \mathbf{R}^d} \sup_{\substack{0 \leq s,t \leq T, \\ 0 < |x'-x| \leq \rho}} \frac{|D_x^\alpha \phi_{s,t}(x,\omega) - D_x^\alpha \phi_{s,t}(x',\omega)|}{|x-x'|^\epsilon (1+|x|)^\gamma},$$

*are finite. Furthermore, the random variables defined by the above expressions have p-th moments for all* $p \geq 1$.

*Proof.*

The cocycle property stated in (ii) is proved in [I-W] for the white-noise case using an approximation argument (cf. [Mo.1], [Mo.2]). Assertions (iii) and (iv) are well-known to hold for a.a. $\omega \in \Omega$ ([Ku], Theorem 4.6.5). A perfect version of $\phi_{s,t}$ satisfying (i)-(iv) for *all* $\omega \in \Omega$, is established in [A-S]. The arguments in [A-S] use perfection techniques and Theorem 4.6.5 of [Ku] (cf. also [M-S.1]).

Assume that for every $\gamma, T > 0$ the two random variables in (v) have finite moments of all orders. Let $\Omega^{T,\gamma}$ be the set of all $\omega \in \Omega$ for which all random variables in (v) are finite. Define the set $\Omega_0$ by

$$\Omega_0 := \bigcap_{T \in \mathbf{N}} \bigcap_{n \in \mathbf{N}} \bigcap_{s \in \mathbf{R}} \theta(s,\cdot)(\Omega^{T,1/n}).$$

Then $\theta(s,\cdot)(\Omega_0) = \Omega_0$ for all $s \in \mathbf{R}$. Furthermore, it is not hard to see that

$$\bigcap_{T \in \mathbf{N}} \bigcap_{n \in \mathbf{N}} \bigcap_{m \in \mathbf{Z}} \theta(mT,\cdot)(\Omega^{2T,1/n}) \subseteq \Omega_0.$$

Therefore $\Omega_0$ is a sure event in $\bar{\mathcal{F}}$. By Proposition 2.3, $\Omega_0$ contains a sure invariant event $\Omega_0' \in \mathcal{F}$. Hence we can redefine $\phi_{s,t}(\cdot,\omega)$ and $\phi(t,\cdot,\omega)$ to be the identity map $\mathbf{R}^d \to \mathbf{R}^d$ for all $\omega \in \Omega \setminus \Omega_0'$. This can be done without violating properties (i)-(iv).

By Proposition 2.2, Theorem 1 in [M-S.2] and the remark following its proof, it follows that the two random variables

$$X_1 := \sup_{\substack{0 \leq s \leq t \leq T, \\ x \in \mathbf{R}^d}} \frac{|\phi_{s,t}(x,\cdot)|}{[1+|x|(\log^+|x|)^\gamma]}, \quad X_2 := \sup_{\substack{0 \leq s \leq t \leq T, \\ x \in \mathbf{R}^d}} \frac{|x|}{[1+|\phi_{s,t}(x,\cdot)|(\log^+|x|)^\gamma]}$$

have $p$-th moments for all $p \geq 1$. To complete the proof of the first assertion in (v), it is sufficient to show that the random variable

$$\hat{X}_1 := \sup_{\substack{0 \leq s \leq t \leq T, \\ x \in \mathbf{R}^d}} \frac{|\phi_{t,s}(x,\cdot)|}{[1+|x|(\log^+|x|)^\gamma]}$$

has $p$-th moments for all $p \geq 1$. To do this, assume (without loss of generality) that $\gamma \in (0,1)$. From the definition of $X_2$, we have

$$|y| \leq X_2[1 + |\phi_{s,t}(y, \cdot)|(\log^+ |y|)^\gamma]$$

for all $0 \leq s \leq t \leq T, y \in \mathbf{R}^d$. Use the substitution

$$y = \phi_{t,s}(x, \omega) = \phi_{s,t}^{-1}(x, \omega), \quad \phi_{s,t}(y, \omega) = x, \quad 0 \leq s \leq t \leq T, \omega \in \Omega, x \in \mathbf{R}^d,$$

to rewrite the above inequality in the form

$$|y| \leq X_2[1 + |x|(\log^+ |y|)^\gamma].$$

By an elementary computation, the above inequality may be solved for $\log^+ |y|$. This gives a positive non-random constant $K_1$ (possibly dependent on $\epsilon$ and $T$) such that

$$|y| \leq K_1 X_2[1 + |x|\{1 + (\log^+ |X_2|)^\gamma + (\log^+ |x|)^\gamma\}].$$

Since $X_2$ has moments of all orders, the above inequality implies that $\hat{X}_1$ also has $p$-th moments for all $p \geq 1$.

We now prove the second assertion in (v). First, note that the following two random variables

$$X_3 := \sup_{\substack{0 \leq s \leq t \leq T, \\ x \in \mathbf{R}^d}} \frac{|D_x^\alpha \phi_{s,t}(x, \cdot)|}{(1 + |x|^\gamma)}, \ |\alpha| \leq k, \quad X_4 := \sup_{\substack{0 \leq s \leq t \leq T, \\ x \in \mathbf{R}^d}} \frac{|[D\phi_{s,t}(x, \cdot)]^{-1}|}{(1 + |x|^\gamma)}$$

have $p$-th moments for all $p \geq 1$ ([Ku], Ex. 4.6.9, p. 176; [M-S.2], Remark (i) following Theorem 2). We must show that the random variables

$$\hat{X}_3 := \sup_{\substack{0 \leq s \leq t \leq T, \\ x \in \mathbf{R}^d}} \frac{|D_x^\alpha \phi_{s,t}^{-1}(x, \cdot)|}{(1 + |x|^\gamma)}, \quad 1 \leq |\alpha| \leq k,$$

have $p$-th moments for all $p \geq 1$. Note that there is a positive constant $C$ such that for any non-singular matrix $A$, one has

$$|(\det A)^{-1}| = |\det(A^{-1})| \leq C\|A^{-1}\|^d.$$

Using this fact and applying Proposition 2.2 with $f := \phi_{s,t}$, $1 \leq s \leq t \leq T$, shows that for every $\delta' > 0$, any $i \in \{1, 2, \cdots, d\}$ and any $1 \leq |\alpha| \leq k$, there exists a random variable $K_{\delta'} \in \bigcap_{p \geq 1} L^p(\Omega, \mathbf{R})$ such that

$$|D_x^\alpha(\phi_{s,t}^{-1})_i(x)| \leq K_{\delta'}(1 + |x|^{\delta'})^{m_{\alpha,i}}.$$

for all $x \in \mathbf{R}^d$ and some positive integer $m_{\alpha,i}$. Now for any given $\epsilon > 0$, choose $\delta' = \dfrac{\gamma}{m_{\alpha,i}}$ to obtain

$$|D_x^\alpha (\phi_{s,t}^{-1})_i(x)| \leq K_{\delta'}(1 + |x|^{\gamma/m_{\alpha,i}})^{m_{\alpha,i}} \leq 2^{m_{\alpha,i}} K_{\delta'}(1 + |x|^\gamma)$$

for all $x \in \mathbf{R}^d$. This shows that $\hat{X}_3$ has $p$-th moments for all $p \geq 1$.

The last estimate in (v) follows from a somewhat lengthy argument. We will only sketch it. First note that for every $p \geq 1$, there exists a constant $c \geq 0$ such that

$$E(|D_x^\alpha \phi_{s,t}(x) - D_x^\alpha \phi_{s',t'}(x')|^{2p}) \leq c(|x - x'|^{2p\delta} + |s - s'|^p + |t - t'|^p)$$

uniformly for all $x, x' \in \mathbf{R}^d$, $0 \leq s \leq t \leq T$ ([Ku], Theorem 4.6.4, pp. 172-173). Using the above estimate, we can employ the inequality of Garsia-Rodemich-Rumsey in its majorising measure version in order to show that the expression

$$\sup_{x \in \mathbf{R}^d} \sup_{\substack{0 \leq s \leq t \leq T, \\ 0 < |x'-x| \leq \rho}} \frac{|D_x^\alpha \phi_{s,t}(x, \omega) - D_x^\alpha \phi_{s,t}(x', \omega)|}{|x - x'|^\epsilon (1 + |x|)^\gamma}$$

has moments of all orders. The argument used to show this is similar to the one used in [I-S]. The application of the Garsia-Rodemich-Rumsey inequality is effected using the following metric on the space $[0, T] \times [0, T] \times \mathbf{R}^d$:

$$d((s, t, x), (s', t', x')) := |x - x'|^\delta + |s - s'|^{1/2} + |t - t'|^{1/2}.$$

Finally, we extend the estimate to cover the sup over *all* $(s, t) \in [0, T] \times [0, T]$ by appealing to Proposition 2.2 and the argument used above to establish the existence of $p$-th moments of $\hat{X}_3$. This completes the proof of the theorem. □

## 3. The local stable manifold theorem.

In this section, we shall maintain the general setting and hypotheses of Section 2.

*Furthermore, we shall assume from now on that the $P$-preserving flow $\theta : \mathbf{R} \times \Omega \to \Omega$ is ergodic.*

For any $\rho > 0$ and $x \in \mathbf{R}^d$ denote by $B(x, \rho)$ the open ball with center $x$ and radius $\rho$ in $\mathbf{R}^d$. Denote by $\bar{B}(x, \rho)$ the corresponding closed ball.

Recall that $(\phi, \theta)$ is the perfect cocycle associated with the trajectories $\phi_{s,t}(x)$ of (S) or (I) (Theorem 2.1).

**Definition 3.1.**

Say that the cocycle $\phi$ has a *stationary trajectory* if there exists an $\mathcal{F}$-measurable random variable $Y : \Omega \to \mathbf{R}^d$ such that

$$\phi(t, Y(\omega), \omega) = Y(\theta(t, \omega)) \tag{1}$$

for all $t \in \mathbf{R}$ and every $\omega \in \Omega$. In the sequel, we will always refer to the stationary trajectory (1) by $\phi(t, Y)$.

If (1) is known to hold on a sure event $\Omega_t$ that may depend on $t$, then there are "perfect" versions of the stationary random variable $Y$ and of the flow $\phi$ such that (1) and the conclusions of Theorem 2.1 hold *for all* $\omega \in \Omega$ (under the hypotheses therein)([Sc]).

We may replace $\omega$ in (1) by $\theta(s, \omega)$, $s \in \mathbf{R}$, to get

$$\phi(t, Y(\theta(s, \omega)), \theta(s, \omega)) = Y(\theta(t+s, \omega)) \tag{2}$$

for all $s, t \in \mathbf{R}$ and every $\omega \in \Omega$.

To illustrate the concept of a stationary trajectory, we give a few simple examples.

**Examples.**
(i) Consider the SDE:

$$d\phi(t) = h(\phi(t)) \, dt + \sum_{i=1}^{m} g_i(\phi(t)) \circ dW_i(t)$$

with vector fields $h, g_i : \mathbf{R}^d \to \mathbf{R}^d, i = 1, \cdots, m$, in $C_b^\infty$ and globally bounded. Suppose $h(x_0) = g_i(x_0) = 0$, $1 \leq i \leq m$ for some fixed $x_0 \in \mathbf{R}^d$. Take $Y(\omega) = x_0$ for all $\omega \in \Omega$. Then $Y$ is a stationary trajectory of the above SDE.

(ii) Consider the affine linear one-dimensional SDE:

$$d\phi(t) = \lambda \phi(t) \, dt + dW(t)$$

where $\lambda > 0$ is fixed, and $W(t) \in \mathbf{R}$ is one-dimensional Brownian motion. Take

$$Y(\omega) := -\int_0^\infty e^{-\lambda u} \, dW(u),$$

$$\theta(t, \omega)(s) = \omega(t+s) - \omega(t).$$

Using integration by parts and variation of parameters, the reader may check that there is a version of $Y$ such that $\phi(t, Y(\omega), \omega) = Y(\theta(t, \omega))$ for all $(t, \omega) \in \mathbf{R} \times \Omega$.

(iii) Consider a 2-dimensional affine linear SDE in $\mathbf{R}^2$:

$$d\phi(t) = A\phi(t) \, dt + G dW(t)$$

with $A$ a fixed hyperbolic $(2 \times 2)$-diagonal matrix

$$A := \begin{pmatrix} \lambda_1 & 0 \\ 0 & \lambda_2 \end{pmatrix}$$

where $\lambda_2 < 0 < \lambda_1$. $G$ is a constant matrix, e.g.

$$G := \begin{pmatrix} g_1 & g_2 \\ g_3 & g_4 \end{pmatrix}$$

where $g_i \in \mathbf{R}, i = 1, 2, 3, 4$. $W := \begin{pmatrix} W_1 \\ W_2 \end{pmatrix}$ is 2-dimensional Brownian motion. Set $Y := \begin{pmatrix} Y_1 \\ Y_2 \end{pmatrix}$ where

$$Y_1 := -g_1 \left[ \int_0^\infty e^{-\lambda_1 u} \, dW_1(u) \right] - g_2 \left[ \int_0^\infty e^{-\lambda_1 u} \, dW_2(u) \right]$$

and

$$Y_2 := g_3 \left[ \int_{-\infty}^0 e^{-\lambda_2 u} \, dW_1(u) \right] + g_4 \left[ \int_{-\infty}^0 e^{-\lambda_2 u} \, dW_2(u) \right].$$

Using variation of parameters and integration by parts (as in (ii)), it is easy to see that $Y$ has a measurable version $Y : \Omega \to \mathbf{R}^2$ which gives a stationary trajectory of the SDE that satisfies (1).

(iv) By Itô's formula, non-linear transforms of the SDE in (iii) under a fixed global diffeomorphism of $R^2$, immediately yield a stationary trajectory of the transformed SDE.

**Lemma 3.1.**

Let the conditions of Theorem 2.1 hold. Assume also that $\log^+ |Y(\cdot)|$ is integrable. Then the cocycle $\phi$ satisfies

$$\int_\Omega \log^+ \sup_{-T \leq t_1, t_2 \leq T} \|\phi(t_2, Y(\theta(t_1, \omega)) + (\cdot), \theta(t_1, \omega))\|_{k,\epsilon} \, dP(\omega) < \infty \qquad (3)$$

for any fixed $0 < T, \rho < \infty$ and any $\epsilon \in (0, \delta)$. The symbol $\|\cdot\|_{k,\epsilon}$ denotes the $C^{k,\epsilon}$-norm on $C^{k,\epsilon}$ mappings $\bar{B}(0, \rho) \to \mathbf{R}^d$. Furthermore, the linearized flow $(D_2\phi(t, Y(\omega), \omega), \theta(t, \omega), t \geq 0)$ is an $L(\mathbf{R}^d)$-valued perfect cocycle and

$$\int_\Omega \log^+ \sup_{-T \leq t_1, t_2 \leq T} \|D_2\phi(t_2, Y(\theta(t_1, \omega)), \theta(t_1, \omega))\|_{L(\mathbf{R}^d)} \, dP(\omega) < \infty \qquad (4)$$

for any fixed $0 < T < \infty$. The forward cocycle $(D_2\phi(t, Y(\omega), \omega), \theta(t, \omega), t > 0)$, has a non-random finite Lyapunov spectrum $\{\lambda_m < \cdots < \lambda_{i+1} < \lambda_i < \cdots < \lambda_2 < \lambda_1\}$. Each Lyapunov exponent $\lambda_i$ has a non-random (finite) multiplicity $q_i$, $1 \leq i \leq m$, and $\sum_{i=1}^{m} q_i = d$. The backward linearized cocycle $(D_2\phi(t, Y(\omega), \omega), \theta(t, \omega), t < 0)$, admits a "backward" non-random finite Lyapunov spectrum defined by

$$\lim_{t \to -\infty} \frac{1}{t} \log |D_2\phi(t, Y(\omega), \omega)(v)|, \quad v \in \mathbf{R}^d,$$

and taking values in $\{-\lambda_i\}_{i=1}^{m}$ with non-random (finite) multiplicities $q_i$, $1 \leq i \leq m$, and $\sum_{i=1}^{m} q_i = d$.

Note that Lemma 3.1 stipulates regularity only on the *forward* characteristics of $\overset{\circ}{F}$ and $F$.

*Proof of Lemma 3.1.*

We first prove (4). Start with the perfect cocycle property for $(\phi, \theta)$:

$$\phi(t_1 + t_2, \cdot, \omega) = \phi(t_2, \cdot, \theta(t_1, \omega)) \circ \phi(t_1, \cdot, \omega), \tag{5}$$

for all $t_1, t_2 \in \mathbf{R}$ and all $\omega \in \Omega$. The perfect cocycle property for $(D_2\phi(t, Y(\omega), \omega), \theta(t, \omega))$ follows directly by taking Fréchet derivatives at $Y(\omega)$ on both sides of (5); viz.

$$\begin{aligned} D_2\phi(t_1 + t_2, Y(\omega), \omega) &= D_2\phi(t_2, \phi(t_1, Y(\omega), \omega), \theta(t_1, \omega)) \circ D_2\phi(t_1, Y(\omega), \omega) \\ &= D_2\phi(t_2, Y(\theta(t_1, \omega)), \theta(t_1, \omega)) \circ D_2\phi(t_1, Y(\omega), \omega) \end{aligned} \tag{6}$$

for all $\omega \in \Omega_0, t_1, t_2 \in \mathbf{R}$. The existence of a fixed discrete spectrum for the linearized cocycle follows the analysis in [Mo.1] and [M-S.1]. This analysis uses the integrability property (4) and the ergodicity of $\theta$. Although (4) is an easy consequence of (6) and Theorem 2.1 (v), it is clear that (3) implies (4). Therefore it is sufficient to establish (3).

In view of (1) and the identity

$$\phi_{t_1, t_1+t_2}(x, \omega) = \phi(t_2, x, \theta(t_1, \omega)), \quad x \in \mathbf{R}^d, t_1, t_2 \in \mathbf{R},$$

(Theorem 2.1(i)), (3) will follow if we show that the following integrals are finite for $0 \leq |\alpha| \leq k$:

$$\int_{\Omega} \log^+ \sup_{\substack{0 \leq s, t \leq T, \\ |x'| \leq \rho}} |D_x^\alpha \phi_{s,t}(\phi_{0,s}(Y(\omega), \omega) + x', \omega)| \, dP(\omega), \tag{7}$$

$$\int_{\Omega} \log^+ \sup_{\substack{0 \leq s, t \leq T, \\ |x, x'| \in \bar{B}(0, \rho), x \neq x'}} \frac{|D_x^\alpha \phi_{s,t}(\phi_{0,s}(Y(\omega), \omega) + x, \omega) - D_x^\alpha \phi_{s,t}(\phi_{0,s}(Y(\omega), \omega) + x', \omega)|}{|x - x'|^\epsilon} \, dP(\omega), \tag{7'}$$

For simplicity of notation, we shall denote random constants by the letters $K_i, i = 1, 2, 3, 4$. Each $K_i, i = 1, 2, 3, 4$, has $p$-th moments for all $p \geq 1$ and may depend on $\rho$ and $T$. The following string of inequalities follows easily from Theorem 2.1 (v).

$$
\begin{aligned}
\log^+ \sup_{\substack{s,t \in [0,T], \\ |x'| \leq \rho}} &|D_x^\alpha \phi_{s,t}(\phi_{0,s}(Y(\omega), \omega) + x', \omega)| \\
&\leq \log^+ \sup_{s \in [0,T]} \left\{ K_1(\omega)[1 + (\rho + |\phi_{0,s}(Y(\omega), \omega)|)^2] \right\} \\
&\leq \log^+ K_2(\omega) + \log^+[1 + 2\rho^2 + K_3(\omega)(1 + |Y(\omega)|^4)] \\
&\leq \log^+ K_4(\omega) + \log[1 + 2\rho^2] + 4\log^+ |Y(\omega)| \qquad (8)
\end{aligned}
$$

for all $\omega \in \Omega$. Now (8) and the integrability hypothesis on $Y$ imply that the integral (7) is finite. The finiteness of $(7')$ follows in a similar manner using Theorem 2.1 (v). This completes the proof of the lemma. □

**Definition 3.2.**

A stationary trajectory $\phi(t, Y)$ of $\phi$ is said to be *hyperbolic* if $E \log^+ |Y(\cdot)| < \infty$ and the linearized cocycle $(D_2\phi(t, Y(\omega), \omega), \theta(t, \omega), t \geq 0)$ has a Lyapunov spectrum $\{\lambda_m < \cdots < \lambda_{i+1} < \lambda_i < \cdots < \lambda_2 < \lambda_1\}$ which *does not contain* 0.

Let $\{U(\omega), S(\omega) : \omega \in \Omega\}$ denote the unstable and stable subspaces for the linearized cocycle $(D_2\phi(t, Y(\cdot), \cdot), \theta(t, \cdot))$ as given by Theorem 5.3 in [M-S.1]. See also [Mo.1] . This requires the integrability property (4).

The following discussion is devoted to the Stratonovich sde (S) and the linearization of the stochastic flow around a stationary trajectory.

**The Linearization.**

In (S), suppose $\overset{\circ}{F}$ is a forward-backward semimartingale helix satisfying Hypotheses $(ST(k,\delta))$ and $(ST^-(k,\delta))$ for some $k \geq 2$ and $\delta \in (0,1]$. Then it follows from Theorem 4.2(i) that the (possibly anticipating) process $\phi(t, Y(\omega), \omega)$ is a *trajectory* of the anticipating Stratonovich sde

$$
\left.\begin{aligned}
d\phi(t, Y) &= \overset{\circ}{F}(\circ dt, \phi(t, Y)), \quad t > 0 \\
\phi(0, Y) &= Y.
\end{aligned}\right\} \qquad (SII)
$$

In the above sde, the Stratonovich differential $\overset{\circ}{F}(\circ dt, \cdot)$ is defined as in Section 4 (Definition 4.1, cf. [Ku], p. 86). The above sde follows immediately by substituting $x = Y(\omega)$ in

$$
\left.\begin{aligned}
d\phi(t, x) &= \overset{\circ}{F}(\circ dt, \phi(t, x)), \quad t > 0 \\
\phi(0) &= x \in \mathbf{R}^d
\end{aligned}\right\} \qquad (SI)
$$

(Theorem 4.2 (i)). This substitution works in spite of the anticipating nature of $\phi(t, Y(\omega), \omega) = Y(\theta(t, \omega))$, because the Stratonovich integral is stable under random *anticipating* substitutions (Theorem 4.1).

Furthermore, we can linearize the sde (S) along the stationary trajectory and then match the solution of the linearized equation with the linearized cocycle $D_2\phi(t, Y(\omega), \omega)$. That is to say, the (possibly non-adapted) process

$$y(t) := D_2\phi(t, Y(\omega), \omega), t \geq 0,$$

satisfies the associated Stratonovich linearized sde

$$\left.\begin{aligned}dy(t) &= D_2\overset{\circ}{F}(\circ dt, Y(\theta(t)))y(t), \quad t > 0 \\ y(0) &= I \in L(\mathbf{R}^d).\end{aligned}\right\} \quad (SIII)$$

In (SIII), the symbol $D_2$ denotes the spatial (Fréchet) derivative of the driving semimartingale along the stationary trajectory $\phi(t, Y(\omega), \omega) = Y(\theta(t, \omega))$ (Theorem 4.2(ii)).

In view of Hypothesis $(ST^-(k, \delta))$ for $k \geq 2, \delta \in (0, 1)$, and Theorem 4.2(iii),(iv), it follows that the backward trajectories $\phi(t, Y)$, $\hat{y}(t) := D_2\phi(t, Y, \cdot)$, $t < 0$, satisfy the backward sde's

$$\left.\begin{aligned}d\phi(t, Y) &= -\overset{\circ}{F}(\circ \hat{dt}, \phi(t, Y)), \quad t < 0 \\ \phi(0, Y) &= Y,\end{aligned}\right\} \quad (SII^-)$$

$$\left.\begin{aligned}d\hat{y}(t) &= -D_2\overset{\circ}{F}(\circ \hat{dt}, \phi(t, Y))\hat{y}(t), \quad t < 0, \\ \hat{y}(0, Y) &= I \in L(\mathbf{R}^d).\end{aligned}\right\} \quad (SIII^-)$$

Note however that the significance of (SIII) is to provide a direct link between the linearized flow $D_2\phi(t, Y(\omega), \omega)$ and the linearized sde. The Stratonovich equation (SII) *does not play a direct role* in the construction of the stable and unstable manifolds (cf. [Wa], Section 4.2). On the other hand, (SII) and $(SII^-)$ provide a dynamic characterization of the stable and unstable manifolds in Theorem 3.1 (a), (d).

In order to apply Ruelle's discrete theorem ([Ru.1], Theorem 5.1, p. 292), we will introduce the following auxiliary cocycle $Z : \mathbf{R} \times \mathbf{R}^d \times \Omega \to \mathbf{R}^d$, which is essentially a "centering" of the flow $\phi$ about the stationary solution:

$$Z(t, x, \omega) := \phi(t, x + Y(\omega), \omega) - Y(\theta(t, \omega)) \tag{9}$$

for $t \in \mathbf{R}, x \in \mathbf{R}^d, \omega \in \Omega$.

**Lemma 3.2.**

Assume the hypotheses of Theorem 2.1. Then $(Z, \theta)$ is a perfect cocycle on $\mathbf{R}^d$ and $Z(t, 0, \omega) = 0$ for all $t \in \mathbf{R}$, and all $\omega \in \Omega$.

*Proof.*

Let $t_1, t_2 \in \mathbf{R}, \omega \in \Omega, x \in \mathbf{R}^d$. Then by the cocycle property for $\phi$ and Definition 3.1, we have

$$\begin{aligned}
Z(t_2, Z(t_1, x, \omega), \theta(t_1, \omega)) &= \phi(t_2, Z(t_1, x, \omega) + Y(\theta(t_1, \omega)), \theta(t_1, \omega)) - Y(\theta(t_2, \theta(t_1, \omega))) \\
&= \phi(t_2, \phi(t_1, x + Y(\omega), \omega), \theta(t_1, \omega)) - Y(\theta(t_2 + t_1, \omega)) \\
&= Z(t_1 + t_2, x, \omega).
\end{aligned}$$

The assertion $Z(t, 0, \omega) = 0$, $t \in \mathbf{R}$, $\omega \in \Omega$, follows directly from the definition of $Z$ and Definition 3.1. □

The next lemma will be needed in order to construct the shift-invariant sure events appearing in the statement of the local stable manifold theorem. The lemma essentially gives "perfect versions" of the ergodic theorem and Kingman's subadditive ergodic theorem.

**Lemma 3.3.**

(i) Let $h : \Omega \to \mathbf{R}^+$ be $\mathcal{F}$-measurable and such that

$$\int_\Omega \sup_{0 \leq u \leq 1} h(\theta(u, \omega)) \, dP(\omega) < \infty.$$

Then there is a sure event $\Omega_1 \in \mathcal{F}$ such that $\theta(t, \cdot)(\Omega_1) = \Omega_1$ for all $t \in \mathbf{R}$, and

$$\lim_{t \to \infty} \frac{1}{t} h(\theta(t, \omega)) = 0$$

for all $\omega \in \Omega_1$.

(ii) Suppose $f : \mathbf{R}^+ \times \Omega \to \mathbf{R} \cup \{-\infty\}$ is a measurable process on $(\Omega, \mathcal{F}, P)$ satisfying the following conditions

(a) $\displaystyle\int_\Omega \sup_{0 \leq u \leq 1} f^+(u, \omega) \, dP(\omega) < \infty, \quad \int_\Omega \sup_{0 \leq u \leq 1} f^+(1 - u, \theta(u, \omega)) \, dP(\omega) < \infty$

(b) $f(t_1 + t_2, \omega) \leq f(t_1, \omega) + f(t_2, \theta(t_1, \omega))$ for all $t_1, t_2 \geq 0$ and **all** $\omega \in \Omega$.

Then there is sure event $\Omega_2 \in \mathcal{F}$ such that $\theta(t, \cdot)(\Omega_2) = \Omega_2$ for all $t \in \mathbf{R}$, and a fixed number $f^* \in \mathbf{R} \cup \{-\infty\}$ such that

$$\lim_{t \to \infty} \frac{1}{t} f(t, \omega) = f^*$$

for all $\omega \in \Omega_2$.

*Proof.*

A proof of (i) is given in [Mo.1], Lemma 5 (iii) with a sure event $\tilde{\Omega}_1 \in \bar{\mathcal{F}}$ such that $\theta(t, \cdot)(\tilde{\Omega}_1) \subseteq \tilde{\Omega}_1$ for all $t \geq 0$. Proposition 2.3 now gives a sure event $\Omega_1 \subseteq \tilde{\Omega}_1$ such that $\Omega_1 \in \mathcal{F}$ and satisfies assertion (i) of the lemma.

Assertion (ii) follow from [Mo.1], Lemma 7, and Proposition 2.3. □

The proof of the local stable-manifold theorem (Theorem 3.1) uses a discretization argument that requires the following lemma.

**Lemma 3.4.**

Assume the hypotheses of Lemma 3.2, and suppose that $\log^+ |Y(\cdot)|$ is integrable. Then there is a sure event $\Omega_3 \in \mathcal{F}$ with the following properties:

(i) $\theta(t, \cdot)(\Omega_3) = \Omega_3$ for all $t \in \mathbf{R}$,

(ii) For every $\omega \in \Omega_3$ and any $x \in \mathbf{R}^d$, the statement

$$\limsup_{n \to \infty} \frac{1}{n} \log |Z(n, x, \omega)| < 0 \tag{10}$$

*implies*

$$\limsup_{t \to \infty} \frac{1}{t} \log |Z(t, x, \omega)| = \limsup_{n \to \infty} \frac{1}{n} \log |Z(n, x, \omega)|. \tag{11}$$

*Proof.*

The integrability condition (3) of Lemma 3.1 implies that

$$\int_\Omega \log^+ \sup_{\substack{0 \leq t_1, t_2 \leq 1, \\ x^* \in \bar{B}(0,1)}} \|D_2 Z(t_1, x^*, \theta(t_2, \omega))\|_{L(\mathbf{R}^d)} \, dP(\omega) < \infty. \tag{12}$$

Therefore by (the perfect version of) the ergodic theorem (Lemma 3.3(i)), there is a sure event $\Omega_3 \in \mathcal{F}$ such that $\theta(t, \cdot)(\Omega_3) = \Omega_3$ for all $t \in \mathbf{R}$, and

$$\lim_{t \to \infty} \frac{1}{t} \log^+ \sup_{\substack{0 \leq u \leq 1, \\ x^* \in \bar{B}(0,1)}} \|D_2 Z(u, x^*, \theta(t, \omega))\|_{L(\mathbf{R}^d)} = 0 \tag{13}$$

for all $\omega \in \Omega_3$.

Let $\omega \in \Omega_3$ and suppose $x \in \mathbf{R}^d$ satisfies (10). Then (10) implies that there exists a positive integer $N_0(x, \omega)$ such that $Z(n, x, \omega) \in \bar{B}(0, 1)$ for all $n \geq N_0$. Let $n \leq t < n+1$ where $n \geq N_0$. Then by the cocycle property for $(Z, \theta)$ and the Mean Value Theorem, we have

$$\sup_{n \leq t \leq n+1} \frac{1}{t} \log |Z(t, x, \omega)|$$
$$\leq \frac{1}{n} \log^+ \sup_{\substack{0 \leq u \leq 1, \\ x^* \in \bar{B}(0,1)}} \|D_2 Z(u, x^*, \theta(n, \omega))\|_{L(\mathbf{R}^d)} + \frac{n}{(n+1)} \frac{1}{n} \log |Z(n, x, \omega)|.$$

Take $\limsup_{n \to \infty}$ in the above relation and use (13) to get

$$\limsup_{t \to \infty} \frac{1}{t} \log |Z(t, x, \omega)| \leq \limsup_{n \to \infty} \frac{1}{n} \log |Z(n, x, \omega)|.$$

The inequality

$$\limsup_{n \to \infty} \frac{1}{n} \log |Z(n, x, \omega)| \leq \limsup_{t \to \infty} \frac{1}{t} \log |Z(t, x, \omega)|,$$

is obvious. Hence (11) holds, and the proof of the lemma is complete. □

In order to formulate the measurability properties of the stable and unstable manifolds, we will consider the class $\mathcal{C}(\mathbf{R}^d)$ of all non-empty compact subsets of $\mathbf{R}^d$. Give $\mathcal{C}(\mathbf{R}^d)$ the Hausdorff metric $d^*$:

$$d^*(A_1, A_2) := \sup\{d(x, A_1) : x \in A_2\} \vee \sup\{d(y, A_2) : y \in A_1\}$$

where $A_1, A_2 \in \mathcal{C}(\mathbf{R}^d)$, and $d(x, A_i) := \inf\{|x - y| : y \in A_i\}$, $x \in \mathbf{R}^d$, $i = 1, 2$. Denote by $\mathcal{B}(\mathcal{C}(\mathbf{R}^d))$ the Borel $\sigma$-algebra on $\mathcal{C}(\mathbf{R}^d)$ with respect to the metric $d^*$. Then $(\mathcal{C}(\mathbf{R}^d), d^*)$ is a complete separable metric space. Morevover, it is not hard to see that finite non-empty intersections are jointly measurable and translations are jointly continuous on $\mathcal{C}(\mathbf{R}^d)$. These facts are used in the proof of Theorem 3.1 (h).

We now state the local stable manifold theorem for the sde's (S) and (I) around a hyperbolic stationary solution.

**Theorem 3.1.** *(Local Stable and Unstable Manifolds)*

*Assume that $\overset{\circ}{F}$ satisfies Hypothesis $(ST(k,\delta))$ (resp. $F$ satisfies $(IT(k,\delta))$) for some $k \geq 1$ and $\delta \in (0, 1]$. Suppose $\phi(t, Y)$ is a hyperbolic stationary trajectory of (S) (resp. (I)) with $E \log^+ |Y| < \infty$. Suppose the linearized cocycle $(D_2\phi(t, Y(\omega), \omega), \theta(t, \omega), t \geq 0)$ has a Lyapunov spectrum $\{\lambda_m < \cdots < \lambda_{i+1} < \lambda_i < \cdots < \lambda_2 < \lambda_1\}$. Define $\lambda_{i_0} := \max\{\lambda_i : \lambda_i < 0\}$ if at least one $\lambda_i < 0$. If all $\lambda_i > 0$, set $\lambda_{i_0} = -\infty$. (This implies that $\lambda_{i_0-1}$ is the smallest positive Lyapunov exponent of the linearized flow, if at least one $\lambda_i > 0$; in case all $\lambda_i$ are negative, set $\lambda_{i_0-1} = \infty$.)*

*Fix $\epsilon_1 \in (0, -\lambda_{i_0})$ and $\epsilon_2 \in (0, \lambda_{i_0-1})$. Then there exist*
(i) *a sure event $\Omega^* \in \mathcal{F}$ with $\theta(t, \cdot)(\Omega^*) = \Omega^*$ for all $t \in \mathbf{R}$,*
(ii) *$\mathcal{F}$-measurable random variables $\rho_i, \beta_i : \Omega^* \to [0, \infty)$, $\beta_i > \rho_i > 0$, $i = 1, 2$, such that for each $\omega \in \Omega^*$, the following is true:*
  *There are $C^{k,\epsilon}$ ($\epsilon \in (0, \delta)$) submanifolds $\tilde{\mathcal{S}}(\omega), \tilde{\mathcal{U}}(\omega)$ of $\bar{B}(Y(\omega), \rho_1(\omega))$ and $\bar{B}(Y(\omega), \rho_2(\omega))$ (resp.) with the following properties:*
(a) *$\tilde{\mathcal{S}}(\omega)$ is the set of all $x \in \bar{B}(Y(\omega), \rho_1(\omega))$ such that*

$$|\phi(n, x, \omega) - Y(\theta(n, \omega))| \leq \beta_1(\omega) e^{(\lambda_{i_0} + \epsilon_1)n}$$

*for all integers $n \geq 0$. Furthermore,*

$$\limsup_{t \to \infty} \frac{1}{t} \log |\phi(t, x, \omega) - Y(\theta(t, \omega))| \leq \lambda_{i_0} \tag{14}$$

*for all $x \in \tilde{\mathcal{S}}(\omega)$. Each stable subspace $\mathcal{S}(\omega)$ of the linearized flow $D_2\phi$ is tangent at $Y(\omega)$ to the submanifold $\tilde{\mathcal{S}}(\omega)$, viz. $T_{Y(\omega)}\tilde{\mathcal{S}}(\omega) = \mathcal{S}(\omega)$. In particular, $\dim \tilde{\mathcal{S}}(\omega) = \dim \mathcal{S}(\omega)$ and is non-random.*

(b) $\limsup\limits_{t\to\infty} \dfrac{1}{t} \log\left[\sup\left\{\dfrac{|\phi(t,x_1,\omega) - \phi(t,x_2,\omega)|}{|x_1 - x_2|} : x_1 \neq x_2,\, x_1, x_2 \in \tilde{\mathcal{S}}(\omega)\right\}\right] \leq \lambda_{i_0}$.

(c) (Cocycle-invariance of the stable manifolds):
There exists $\tau_1(\omega) \geq 0$ such that
$$\phi(t,\cdot,\omega)(\tilde{\mathcal{S}}(\omega)) \subseteq \tilde{\mathcal{S}}(\theta(t,\omega)), \quad t \geq \tau_1(\omega). \tag{15}$$
Also
$$D_2\phi(t, Y(\omega), \omega)(\mathcal{S}(\omega)) = \mathcal{S}(\theta(t,\omega)), \quad t \geq 0. \tag{16}$$

(d) $\tilde{\mathcal{U}}(\omega)$ is the set of all $x \in \bar{B}(Y(\omega), \rho_2(\omega))$ with the property that
$$|\phi(-n, x, \omega) - Y(\theta(-n, \omega))| \leq \beta_2(\omega)\, e^{(-\lambda_{i_0-1}+\epsilon_2)n} \tag{17}$$
for all integers $n \geq 0$. Also
$$\limsup\limits_{t\to\infty} \dfrac{1}{t} \log |\phi(-t, x, \omega) - Y(\theta(-t, \omega))| \leq -\lambda_{i_0-1}. \tag{18}$$
for all $x \in \tilde{\mathcal{U}}(\omega)$. Furthermore, $\mathcal{U}(\omega)$ is the tangent space to $\tilde{\mathcal{U}}(\omega)$ at $Y(\omega)$. In particular, $\dim \tilde{\mathcal{U}}(\omega) = \dim \mathcal{U}(\omega)$ and is non-random.

(e) $\limsup\limits_{t\to\infty} \dfrac{1}{t} \log\left[\sup\left\{\dfrac{|\phi(-t,x_1,\omega) - \phi(-t,x_2,\omega)|}{|x_1 - x_2|} : x_1 \neq x_2,\, x_1, x_2 \in \tilde{\mathcal{U}}(\omega)\right\}\right] \leq -\lambda_{i_0-1}$.

(f) (Cocycle-invariance of the unstable manifolds):
There exists $\tau_2(\omega) \geq 0$ such that
$$\phi(-t,\cdot,\omega)(\tilde{\mathcal{U}}(\omega)) \subseteq \tilde{\mathcal{U}}(\theta(-t,\omega)), \quad t \geq \tau_2(\omega). \tag{19}$$
Also
$$D_2\phi(-t, Y(\omega), \omega)(\mathcal{U}(\omega)) = \mathcal{U}(\theta(-t,\omega)), \quad t \geq 0. \tag{20}$$

(g) The submanifolds $\tilde{\mathcal{U}}(\omega)$ and $\tilde{\mathcal{S}}(\omega)$ are transversal, viz.
$$\mathbf{R}^d = T_{Y(\omega)}\tilde{\mathcal{U}}(\omega) \oplus T_{Y(\omega)}\tilde{\mathcal{S}}(\omega). \tag{21}$$

(h) The mappings
$$\Omega \to \mathcal{C}(\mathbf{R}^d), \qquad \Omega \to \mathcal{C}(\mathbf{R}^d),$$
$$\omega \mapsto \tilde{\mathcal{S}}(\omega) \qquad \omega \mapsto \tilde{\mathcal{U}}(\omega)$$
are $(\mathcal{F}, \mathcal{B}(\mathcal{C}(\mathbf{R}^d)))$-measurable.

Assume, in addition, that $\overset{\circ}{F}$ satisfies Hypothesis $(ST(k,\delta))$ (resp. $F$ satisfies $(IT(k,\delta))$) for every $k \geq 1$ and $\delta \in (0,1]$. Then the local stable and unstable manifolds $\tilde{\mathcal{S}}(\omega), \tilde{\mathcal{U}}(\omega)$ are $C^\infty$.

The following corollary follows from Theorem 3.1. See [Ru.1], Section (5.3), p. 49.

**Corollary 3.1.1.** *(White-noise, Itô case)*
  *Consider the Itô sde*

$$dx(t) = h(x(t))\, dt + \sum_{i=1}^{m} g_i(x(t)) dW_i(t) \qquad (V)$$

*Suppose that for some $k \geq 1, \delta \in (0,1)$, $h, g_i$, $1 \leq i \leq m$, are $C_b^{k,\delta}$ vector fields on $\mathbf{R}^d$, and $W := (W_1, \cdots, W_m)$ is m-dimensional Brownian motion on Wiener space $(\Omega, \mathcal{F}, P)$. Let $\theta : \mathbf{R} \times \Omega \to \Omega$ denote the canonical Brownian shift*

$$\theta(t, \omega)(s) := \omega(t+s) - \omega(t), \quad t, s \in \mathbf{R}, \omega \in \Omega. \qquad (22)$$

*Suppose $\phi(t, Y)$ is a hyperbolic stationary trajectory of (V) with $E \log^+ |Y| < \infty$. Then the conclusions of Theorem 3.1 hold.*

  *Furthermore, if the vector fields $h, g_i$, $1 \leq i \leq m$, are $C_b^\infty$, then the conclusions of Theorem 3.1 hold, where $\tilde{\mathcal{S}}(\omega), \tilde{\mathcal{U}}(\omega)$ are $C^\infty$ manifolds.*

**Remarks.**
(i) A similar statement to that of Corollary 3.1.1 holds for the corresponding Stratonovich sde driven by finite-dimensional Brownian motion:

$$dx(t) = h(x(t))\, dt + \sum_{i=1}^{m} g_i(x(t)) \circ dW_i(t) \qquad (SIV)$$

However, in this case one needs stronger conditions to ensure that Hypothesis $(ST(k, \delta))$ hold for (SIV). In fact, such hypotheses will hold if we assume that the functions

$$\mathbf{R}^d \ni x \mapsto \sum_{l=1}^{m} \frac{\partial g_l^i(x)}{\partial x_j} g_l^j(x) \in \mathbf{R}$$

are in $C_b^{k,\delta}$ for each $1 \leq i, j \leq d$ and some $k \geq 1, \delta \in (0,1)$. Cf. the conditions expressed in [A-I]. For example this holds if for some $k \geq 1, \delta \in (0,1)$, the vector field $h$ is of class $C_b^{k,\delta}$ and $g_i$, $1 \leq i \leq m$, are globally bounded and of class $C_b^{k+1,\delta}$. We conjecture that the the global boundedness condition is not needed. This conjecture is not hard to check if the vector fields $g_i$, $1 \leq i \leq m$ are $C_b^\infty$ and generate a finite-dimensional solvable Lie algebra. See [Ku], Theorem 4.9.10, p. 212.
(ii) Recall that if $\overset{\circ}{F}$ is a forward-backward semimartingale helix satisfying Hypotheses $(ST(k, \delta))$ and $(ST^-(k, \delta))$ for some $k \geq 2$ and $\delta \in (0, 1)$, then the inverse $\phi(t, \cdot, \theta(-t, \omega))^{-1}(x) = \phi(-t, x, \omega)$, $t > 0$, corresponds to a *solution* of the sde $(S^-)$. Furthermore, $\phi(-t, Y)$ and $D_2\phi(-t, Y)$, $t > 0$, satisfy the anticipating sde's $(SII^-)$ and $(SIII^-)$, respectively. See Theorem 4.2 (iii), (iv), of Section 4.

(iii) We may replace the stationary random variable $Y$ by its invariant distribution $\mu$ and then formulate all our results with respect to the product measure $\mu \otimes P$ and the underlying skew-product flow. This would give stable and unstable manifolds that are defined a.e.($\mu \otimes P$); cf. [C] for the globally asymptotically stable case on a compact manifold.

(iv) In Corollary 3.1.1, one can allow for infinitely many Brownian motions (cf. [Ku], p. 106-107). Details are left to the reader.

*Proof of Theorem 3.1.*

Assume the hypotheses of the theorem. Consider the cocycle $(Z, \theta)$ defined by (9). Define the family of maps $F_\omega : \mathbf{R}^d \to \mathbf{R}^d$ by $F_\omega(x) := Z(1, x, \omega)$ for all $\omega \in \Omega$ and $x \in \mathbf{R}^d$. Let $\tau := \theta(1, \cdot) : \Omega \to \Omega$. Following Ruelle ([Ru.1], p. 292), define $F_\omega^n := F_{\tau^{n-1}(\omega)} \circ \cdots \circ F_{\tau(\omega)} \circ F_\omega$. Then by the cocycle property for $Z$, we get $F_\omega^n = Z(n, \cdot, \omega)$ for each $n \geq 1$. Clearly, each $F_\omega$ is $C^{k,\epsilon}$ ($\epsilon \in (0, \delta)$) and $(DF_\omega)(0) = D_2\phi(1, Y(\omega), \omega)$. By measurability of the flow $\phi$, it follows that the map $\omega \mapsto (DF_\omega)(0)$ is $\mathcal{F}$-measurable. By (4) of Lemma 3.1, it is clear that the map $\omega \mapsto \log^+ \|D_2\phi(1, Y(\omega), \omega)\|_{L(\mathbf{R}^d)}$ is integrable. Furthermore, the discrete cocycle $((DF_\omega^n)(0), \theta(n, \omega), n \geq 0)$ has a non-random Lyapunov spectrum which coincides with that of the linearized continuous cocycle $(D_2\phi(t, Y(\omega), \omega), \theta(t, \omega), t \geq 0)$, viz. $\{\lambda_m < \cdots < \lambda_{i+1} < \lambda_i < \cdots < \lambda_2 < \lambda_1\}$, where each $\lambda_i$ has fixed multiplicity $q_i$, $1 \leq i \leq m$ (Lemma 3.1). Note that $\lambda_{i_0}$ (and $\lambda_{i_0-1}$) are well-defined by hyperbolicity of the stationary trajectory. If $\lambda_i > 0$ for all $1 \leq i \leq m$, then take $\tilde{\mathcal{S}}(\omega) := \{Y(\omega)\}$ for all $\omega \in \Omega$. The assertions of the theorem are trivial in this case. From now on suppose that at least one $\lambda_i < 0$.

We use Theorem 5.1 of Ruelle ([Ru.1], p. 292) and its proof to obtain a sure event $\Omega_1^* \in \mathcal{F}$ such that $\theta(t, \cdot)(\Omega_1^*) = \Omega_1^*$ for all $t \in \mathbf{R}$, $\mathcal{F}$-measurable positive random variables $\rho_1, \beta_1 : \Omega_1^* \to (0, \infty)$, $\rho_1 < \beta_1$, and a random family of $C^{k,\epsilon}$ ($\epsilon \in (0, \delta)$) submanifolds of $\bar{B}(0, \rho_1(\omega))$ denoted by $\tilde{\mathcal{S}}_d(\omega)$, $\omega \in \Omega_1^*$, and satisfying the following properties for each $\omega \in \Omega_1^*$:

$$\tilde{\mathcal{S}}_d(\omega) = \{x \in \bar{B}(0, \rho_1(\omega)) : |Z(n, x, \omega)| \leq \beta_1(\omega) e^{(\lambda_{i_0} + \epsilon_1)n} \text{ for all integers } n \geq 0\}. \tag{23}$$

Each $\tilde{\mathcal{S}}_d(\omega)$ is tangent at 0 to the stable subspace $\mathcal{S}(\omega)$ of the linearized flow $D_2\phi$, viz. $T_0\tilde{\mathcal{S}}_d(\omega) = \mathcal{S}(\omega)$. In particular, $\dim \tilde{\mathcal{S}}_d(\omega)$ is non-random by the ergodicity of $\theta$. Furthermore,

$$\limsup_{n \to \infty} \frac{1}{n} \log \left[ \sup_{\substack{x_1 \neq x_2, \\ x_1, x_2 \in \tilde{\mathcal{S}}_d(\omega)}} \frac{|Z(n, x_1, \omega) - Z(n, x_2, \omega)|}{|x_1 - x_2|} \right] \leq \lambda_{i_0}. \tag{24}$$

Before we proceed with the proof, we will indicate how one may arrive at the above $\theta(t, \cdot)$-invariant sure event $\Omega_1^* \in \mathcal{F}$ from Ruelle's proof. Consider the proof of Theorem 5.1 in [Ru.1], p. 293. In the notation of [Ru.1], set $T_\omega^t := D_2Z(t, 0, \omega), T_n(\omega) := D_2Z(1, 0, \theta(n-1, \omega)), \tau^t(\omega) := \theta(t, \omega)$, for $t \in \mathbf{R}^+, n = 1, 2, 3, \cdots$. By the integrability condition (4) (Lemma 3.1) and Lemma 3.3 (i),(ii), there is a sure event $\Omega_1^* \in \mathcal{F}$ such that $\theta(t, \cdot)(\Omega_1^*) = \Omega_1^*$ for all $t \in \mathbf{R}$, with the property that *continuous-time* analogues of equations (5.2), (5.3), (5.4) in ([Ru.1], p. 45) hold. In particular,

$$\left. \begin{aligned} \lim_{t \to \infty} \{[D_2 Z(t, 0, \omega)]^* [D_2 Z(t, 0, \omega)]\}^{1/(2t)} &= \Lambda(\omega), \\ \lim_{t \to \infty} \frac{1}{t} \log^+ \|Z(1, \cdot, \theta(t, \omega))\|_{1,\epsilon} &= 0, \end{aligned} \right\} \quad (25)$$

for all $\omega \in \Omega_1^*$, $\epsilon \in (0, \delta)$. See Theorem (B.3) ([Ru.1], p. 304). The rest of the proof of Theorem 5.1 works for a *fixed choice* of $\omega \in \Omega_1^*$. In particular, (the proof of) the "perturbation theorem" ([Ru.1], Theorem 4.1) does not affect the choice of the sure event $\Omega_1^*$ because it works *pointwise* in $\omega \in \Omega_1^*$, and hence does not involve the selection of a sure event ([Ru.1], pp. 285-292).

For each $\omega \in \Omega_1^*$, let $\tilde{\mathcal{S}}(\omega)$ be the set defined in part (a) of the theorem. Then it is easy to see from (23) and the definition of $Z$ that

$$\tilde{\mathcal{S}}(\omega) = \tilde{\mathcal{S}}_d(\omega) + Y(\omega). \quad (26)$$

Since $\tilde{\mathcal{S}}_d(\omega)$ is a $C^{k,\epsilon}$ ($\epsilon \in (0, \delta)$) submanifold of $\bar{B}(0, \rho_1(\omega))$, it follows from (26) that $\tilde{\mathcal{S}}(\omega)$ is a $C^{k,\epsilon}$ ($\epsilon \in (0, \delta)$) submanifold of $\bar{B}(Y(\omega), \rho_1(\omega))$. Furthermore, $T_{Y(\omega)} \tilde{\mathcal{S}}(\omega) = T_0 \tilde{\mathcal{S}}_d(\omega) = \mathcal{S}(\omega)$. In particular, $\dim \tilde{\mathcal{S}}(\omega) = \dim \mathcal{S}(\omega) = \sum_{i=i_0}^{m} q_i$, and is non-random.

Now (24) implies that

$$\limsup_{n \to \infty} \frac{1}{n} \log |Z(n, x, \omega)| \leq \lambda_{i_0} \quad (27)$$

for all $\omega$ in the shift-invariant sure event $\Omega_1^*$ and all $x \in \tilde{\mathcal{S}}_d(\omega)$. Therefore by Lemma 3.4, there is a sure event $\Omega_2^* \subseteq \Omega_1^*$ such that $\theta(t, \cdot)(\Omega_2^*) = \Omega_2^*$ for all $t \in \mathbf{R}$, and

$$\limsup_{t \to \infty} \frac{1}{t} \log |Z(t, x, \omega)| \leq \lambda_{i_0} \quad (28)$$

for all $\omega \in \Omega_2^*$ and all $x \in \tilde{\mathcal{S}}_d(\omega)$. This immediately implies assertion (14) of the theorem.

To prove assertion (b) of the theorem, let $\omega \in \Omega_1^*$. By (24), there is a positive integer $N_0 := N_0(\omega)$ (independent of $x \in \tilde{\mathcal{S}}_d(\omega)$) such that $Z(n, x, \omega) \in \bar{B}(0, 1)$ for all $n \geq N_0$. Let $\Omega_4^* := \Omega_2^* \cap \Omega_3$, where $\Omega_3$ is the shift-invariant sure event defined in the proof of Lemma 3.4. Then $\Omega_4^*$ is a sure event and $\theta(t, \cdot)(\Omega_4^*) = \Omega_4^*$ for all

$t \in \mathbf{R}$. Using an argument similar to the one used in the proof of Lemma 3.4, it follows that

$$\sup_{n \leq t \leq n+1} \frac{1}{t} \log \left[ \sup_{\substack{x_1 \neq x_2, \\ x_1, x_2 \in \tilde{\mathcal{S}}(\omega)}} \frac{|\phi(t, x_1, \omega) - \phi(t, x_2, \omega)|}{|x_1 - x_2|} \right]$$

$$= \sup_{n \leq t \leq n+1} \frac{1}{t} \log \left[ \sup_{\substack{x_1 \neq x_2, \\ x_1, x_2 \in \tilde{\mathcal{S}}_d(\omega)}} \frac{|Z(t, x_1, \omega) - Z(t, x_2, \omega)|}{|x_1 - x_2|} \right]$$

$$\leq \frac{1}{n} \log^+ \sup_{\substack{0 \leq u \leq 1, \\ x^* \in \bar{B}(0,1)}} \|D_2 Z(u, x^*, \theta(n, \omega))\|_{L(\mathbf{R}^d)}$$

$$+ \frac{n}{(n+1)} \frac{1}{n} \log \left[ \sup_{\substack{x_1 \neq x_2, \\ x_1, x_2 \in \tilde{\mathcal{S}}_d(\omega)}} \frac{|Z(n, x_1, \omega) - Z(n, x_2, \omega)|}{|x_1 - x_2|} \right]$$

for all $\omega \in \Omega_4^*$, all $n \geq N_0(\omega)$ and sufficiently large. Taking $\limsup_{n \to \infty}$ in the above inequality and using (24), immediately gives assertion (b) of the theorem.

To prove the invariance property (16), we apply the Oseledec theorem to the linearized cocycle $(D_2\phi(t, Y(\omega), \omega), \theta(t, \omega))$ ([Mo.1], Theorem 4, Corollary 2). This gives a sure $\theta(t, \cdot)$-invariant event, also denoted by $\Omega_1^*$, such that

$$D_2\phi(t, Y(\omega), \omega)(\mathcal{S}(\omega)) \subseteq \mathcal{S}(\theta(t, \omega))$$

for all $t \geq 0$ and all $\omega \in \Omega_1^*$. Equality holds because $D_2\phi(t, Y(\omega), \omega)$ is injective and $\dim \mathcal{S}(\omega) = \dim \mathcal{S}(\theta(t, \omega))$ for all $t \geq 0$ and all $\omega \in \Omega_1^*$.

To prove the asymptotic invariance property (15), we will need to take a closer look at the proofs of Theorems 5.1 and 4.1 in [Ru.1], pp. 285-297. We will first show that $\rho_1, \beta_1$ and a sure event (also denoted by) $\Omega_1^*$ may be chosen such that $\theta(t, \cdot)(\Omega_1^*) = \Omega_1^*$ for all $t \in \mathbf{R}$, and with the property that for any $\epsilon \in (0, \epsilon_1)$ and every $\omega \in \Omega_1^*$, there exists a positive $K_1^\epsilon(\omega)$ for which the inequalities

$$\rho_1(\theta(t, \omega)) \geq K_1^\epsilon(\omega) \rho_1(\omega) e^{(\lambda_{i_0} + \epsilon)t}, \quad \beta_1(\theta(t, \omega)) \geq K_1^\epsilon(\omega) \beta_1(\omega) e^{(\lambda_{i_0} + \epsilon)t} \qquad (29)$$

hold for all $t \geq 0$. The above inequalities hold in the *discrete* case (when $t = n$, a positive integer) from Ruelle's theorem ([Ru.1], Remark (c), p. 297, following the proof of Theorem 5.1). We claim that the relations (29) hold also for *continuous* time. To see this, we will use the method of proof of Theorems 5.1 and 4.1 in [Ru.1]. In the notation of the proof of Theorem 5.1 ([Ru.1], p. 293), observe that the random variable $G$ in (5.5) may be replaced by the larger one

$$\tilde{G}(\omega) := \sup_{t \geq 0} \|F_{\tau^t \omega}\|_{1,\theta} \, e^{(-t\eta - \lambda\theta)} < +\infty, \quad \theta \in (0, 1] \qquad (30)$$

for $0 < \eta < -(\lambda_{i_0} + \epsilon)/4$, and Ruelle's $\lambda$ corresponds to $\lambda_{i_0} + \epsilon_1$ in our notation. Now $\beta_1$ may be chosen using $\delta, A$ from Theorem (4.1) ([Ru.1]) and replacing $G$ with $\tilde{G}$ in (5.10) ([Ru.1], p. 293). Note that $\rho_1 := \frac{\beta_1}{B'_\epsilon}$, where $B'_\epsilon$ is given by (4.5) in Theorem (4.1) of ([Ru.1], p. 285). Therefore, given any *fixed* $\omega \in \Omega_1^*$, we need to determine how the choices of Ruelle's constants $\delta$, $A$ and $B'_\epsilon$ are affected if $\omega$ is replaced by $\theta(l,\omega) = \tau^l(\omega)$ where $l$ is any positive *real* number. Since $T_n(\tau^l(\omega)) = D_2 Z(1, 0, \theta(n-1, \theta(l, \omega))) = T_{n+l}(\omega)$, for all positive integers $n$, it is sufficient to apply Theorem 4.1 ([Ru.1]) to the sequence $\{T_{n+l}(\omega)\}_{n=1}^\infty$. Hence we may follow the discussion in Section (4.7) ([Ru.1], pp. 291-292). We claim that the argument therein still works for positive real $l$. We will indicate the reasoning for $\delta$ and leave the rest of the details to the reader. Consider the definition of $\delta$ in (4.15) in the proof of Lemma (4.2) ([Ru.1], p. 288). Set $\delta(\omega) := \delta, D(\omega) := D, C(\omega) := C$ given by (4.15), (4.11), (4.13), respectively. Redefine $D$ and $C$ by larger constants which we will denote by the same symbols:

$$D(\omega) := \sup_{t \geq 0} e^{-t\eta} \|(\xi^{(t)})^{-1}\| < \infty, \tag{31}$$

$$C(\omega) := \sup_{\substack{0 < s < t < \infty \\ 1 \leq h, k \leq m}} \left\{ \frac{\|T^t \xi_h^{(0)}\| \|T^s \xi_k^{(0)}\|}{\|T^s \xi_h^{(0)}\| \|T^t \xi_k^{(0)}\|} \exp\{[\lambda^{(r(k))} - \lambda^{(r(h))}](t-s)\} \right\} < \infty, \tag{32}$$

where $\xi^{(t)} := (\xi_1^{(t)}, \cdots, \xi_m^{(t)})$, $\xi_k^{(t)} := \frac{T^t \xi_k^{(0)}}{\|T^t \xi_k^{(0)}\|}$, $1 \leq k \leq m$. The $\lambda^{(r(k))}$ are the eigenvalues of $\log \Lambda(\omega)$ with multiplicities. Observe that $D(\omega)$ is finite because the following continuous-time version of (4.9)([Ru.1], p. 287)

$$\lim_{t \to \infty} \frac{1}{t} \log \|T^t \xi_1^{(0)} \wedge \cdots \wedge T^t \xi_m^{(0)}\| = \sum_{k=1}^m \lambda^{(r(k))} \tag{33}$$

holds everywhere on a $\theta(t, \cdot)$-invariant sure event in $\mathcal{F}$ also denoted by $\Omega_1^*$. This is an immediate consequence of Lemma 3.3 (ii). Cf. [Ru.1], p. 287 and p. 303. The constant $C(\omega)$ satisfies the inequality (4.13) of [Ru.1], p. 288, because (by choice of $(t_k^{(n)})^* := t_k^{(n)} e^{[\lambda^{(r(\mu))} - \lambda^{(r(k))}]}$) one has

$$\prod_{k=1}^N (t_k^{(n)})^* = \|T^N \xi_k^{(0)}\| e^{N[\lambda^{(r(\mu))} - \lambda^{(r(k))}]},$$

for all positive integers $N$ and $1 \leq k \leq \mu \leq m$. Indeed, $C(\omega)$ is finite because

$$\lim_{t \to \infty} \frac{1}{t} \log \|T^t \xi_k^{(0)}\| = \lambda^{(r(k))}$$

on a $\theta(t, \cdot)$-invariant sure event in $\mathcal{F}$ (also denoted by $\Omega_1^*$). Now replace $\omega$ in (31) and (32) by $\theta(l, \omega)$. This changes $\xi^{(t)}$ to $\xi^{(t+l)}$, and $T^t \xi_h^{(0)}$ to $T^{t+l} \xi_h^{(0)}$. Hence we get positive constants $K_2^\epsilon(\omega), K_3^\epsilon(\omega)$ such that

$$D(\theta(l,\omega)) \leq K_2^\epsilon(\omega) e^{l(\lambda_{i_0} + \epsilon)} D(\omega), \quad C(\theta(l,\omega)) \leq K_2^\epsilon(\omega) e^{l(\lambda_{i_0} + \epsilon)} C(\omega)$$

STABLE MANIFOLD THEOREM FOR SDE'S 25for sufficiently small $\epsilon \in (0, \epsilon_1)$ and all sufficiently large $l$. From (4.5) in ([Ru.1], p. 288), we get a positive constant $K_4^\epsilon(\omega)$ such that

$$\delta(\theta(l, \omega)) \geq K_4^\epsilon(\omega) e^{l(\lambda_{i_0} + \epsilon)} \delta(\omega)$$

for all sufficiently large $l \in \mathbf{R}^+$ and all sufficiently small $\epsilon$. The behavior of the constants $A$ and $B_\epsilon'$ in Theorem 4.1 ([Ru.1], p. 285) can be analysed in a similar fashion. See [Ru.1], Section (4.7). This yields the inequalities (29). We now proceed to prove (15). Use (b) to obtain a sure event $\Omega_5^* \subseteq \Omega_4^*$ such that $\theta(t, \cdot)(\Omega_5^*) = \Omega_5^*$ for all $t \in \mathbf{R}$, and for any $0 < \epsilon < \epsilon_1$ and $\omega \in \Omega_4^*$, there exists $\beta^\epsilon(\omega) > 0$ (independent of $x$) with

$$|\phi(t, x, \omega) - Y(\theta(t, \omega))| \leq \beta^\epsilon(\omega) e^{(\lambda_{i_0} + \epsilon)t} \tag{34}$$

for all $x \in \tilde{\mathcal{S}}(\omega)$, $t \geq 0$. Fix any real $t \geq 0$, $\omega \in \Omega_5^*$ and $x \in \tilde{\mathcal{S}}(\omega)$. Let $n$ be a non-negative integer. Then the cocycle property and (34) imply that

$$|\phi(n, \phi(t, x, \omega), \theta(t, \omega)) - Y(\theta(n, \theta(t, \omega)))| = |\phi(n+t, x, \omega) - Y(\theta(n+t, \omega))|$$
$$\leq \beta^\epsilon(\omega) e^{(\lambda_{i_0} + \epsilon)(n+t)}$$
$$\leq \beta^\epsilon(\omega) e^{(\lambda_{i_0} + \epsilon)t} e^{(\lambda_{i_0} + \epsilon_1)n}. \tag{35}$$

If $\omega \in \Omega_5^*$, then it follows from (29),(34), (35) and the definition of $\tilde{\mathcal{S}}(\theta(t, \omega))$ that there exists $\tau_1(\omega) > 0$ such that $\phi(t, x, \omega) \in \tilde{\mathcal{S}}(\theta(t, \omega))$ for all $t \geq \tau_1(\omega)$. This proves (15) and completes the proof of assertion (c) of the theorem.

Note that assertions (a), (b) and (c) still hold for all $\omega \in \Omega_5^*$.

We now prove assertion (d) of the theorem, regarding the existence of the local unstable manifolds $\tilde{\mathcal{U}}(\omega)$. We do this by running both the flow $\phi$ and the shift $\theta$ backward in time. Define

$$\tilde{\phi}(t, x, \omega) := \phi(-t, x, \omega), \quad \tilde{Z}(t, x, \omega) := Z(-t, x, \omega), \quad \tilde{\theta}(t, \omega) := \theta(-t, \omega)$$

for all $t \geq 0$ and all $\omega \in \Omega$. Clearly $(\tilde{Z}(t, \cdot, \omega), \tilde{\theta}(t, \omega), t \geq 0)$ is a smooth cocycle, with $\tilde{Z}(t, 0, \omega) = 0$ for all $t \geq 0$. By the hypothesis on $F$ and $Y$, it follows that the linearized flow $(D_2\tilde{\phi}(t, Y(\omega), \omega), \tilde{\theta}(t, \omega), t \geq 0)$ is an $L(\mathbf{R}^d)$-valued perfect cocycle with a non-random finite Lyapunov spectrum $\{-\lambda_1 < -\lambda_2 < \cdots < -\lambda_i < -\lambda_{i+1} < \cdots < -\lambda_m\}$ where $\{\lambda_m < \cdots < \lambda_{i+1} < \lambda_i < \cdots < \lambda_2 < \lambda_1\}$ is the Lyapunov spectrum of the forward linearized flow $(D_2\phi(t, Y(\omega), \omega), \theta(t, \omega), t \geq 0)$. Now apply the first part of the proof of this theorem. This gives *stable manifolds* for the backward flow $\tilde{\phi}$ satisfying assertions (a), (b), (c). This immediately translates into the existence of *unstable manifolds* for the original flow $\phi$, and assertions (d), (e), (f) automatically hold. In particular, we get a sure event $\Omega_6^* \in \mathcal{F}$ such that $\theta(-t, \cdot)(\Omega_6^*) = \Omega_6^*$ for all $t \in \mathbf{R}$, and with the property that assertions (d), (e) and (f) hold for all $\omega \in \Omega_6^*$.

Define the sure event $\Omega^* := \Omega_6^* \cap \Omega_5^*$. Then $\theta(t, \cdot)(\Omega^*) = \Omega^*$ for all $t \in \mathbf{R}$. Furthermore, assertions (a)-(f) hold for all $\omega \in \Omega^*$.

Assertion (g) follows directly from the following facts

$$T_{Y(\omega)}\tilde{\mathcal{U}}(\omega) = \mathcal{U}(\omega), \quad T_{Y(\omega)}\tilde{\mathcal{S}}(\omega) = \mathcal{S}(\omega), \quad \mathbf{R}^d = \mathcal{U}(\omega) \oplus \mathcal{S}(\omega)$$

for all $\omega \in \Omega^*$.

We shall now prove assertion (h). Recall that by (26),

$$\tilde{\mathcal{S}}(\omega) = \mathcal{T}(Y(\omega), \tilde{\mathcal{S}}_d(\omega)) \tag{36}$$

for all $\omega \in \Omega_1^*$, where $\mathcal{T} : \mathbf{R}^d \times \mathcal{C}(\mathbf{R}^d) \to \mathcal{C}(\mathbf{R}^d)$ denotes the translation map

$$\mathcal{T}(x, A) := x + A, \quad x \in \mathbf{R}^d, \ A \in \mathcal{C}(\mathbf{R}^d).$$

Hence, by joint continuity of $\mathcal{T}$ and measurability of $Y$, the $\mathcal{F}$-measurability of the mapping $\Omega \ni \omega \mapsto \tilde{\mathcal{S}}(\omega) \in \mathcal{C}(\mathbf{R}^d)$ would follow from (36) if we can show that the map $\Omega \ni \omega \mapsto \tilde{\mathcal{S}}_d(\omega) \in \mathcal{C}(\mathbf{R}^d)$ is $\mathcal{F}$-measurable. The rest of the argument will demonstrate this.

Define the sequence of random diffeomorphisms

$$f_n(x, \omega) := \beta_1(\omega)^{-1} e^{-(\lambda_{i_0} + \epsilon_1)n} Z(n, x, \omega), \quad x \in \mathbf{R}^d, \ \omega \in \Omega_1^*,$$

for all integers $n \geq 0$. Let $\text{Hom}(\mathbf{R}^d)$ be the topological group of all homeomorphisms of $\mathbf{R}^d$ onto itself. $\text{Hom}(\mathbf{R}^d)$ carries the topology of uniform convergence of sequences of maps and their inverses on compacta. The joint measurability of $f_n$ implies that for each positive integer $n$, the map $\Omega \ni \omega \mapsto f_n(\cdot, \omega) \in \text{Hom}(\mathbf{R}^d)$ is measurable into the Borel field of $\text{Hom}(\mathbf{R}^d)$. Using (23), $\tilde{\mathcal{S}}_d(\omega)$ can be expressed in the form

$$\tilde{\mathcal{S}}_d(\omega) = \lim_{m \to \infty} \bar{B}(0, \rho_1(\omega)) \cap \bigcap_{i=1}^m f_i(\cdot, \omega)^{-1}(\bar{B}(0, 1)) \tag{37}$$

for all $\omega \in \Omega_1^*$. In (37), the limit is taken in the metric $d^*$ on $\mathcal{C}(\mathbf{R}^d)$. The $\mathcal{F}$-measurability of the map $\omega \mapsto \tilde{\mathcal{S}}_d(\omega)$ follows directly from (37), the measurability of $f_i$, $\rho_1$, that of finite intersections and the continuity of the maps

$$\mathbf{R}^+ \ni r \mapsto \bar{B}(0, r) \in \mathcal{C}(\mathbf{R}^d).$$

$$\text{Hom}(\mathbf{R}^d) \ni f \mapsto f^{-1}(\bar{B}(0, 1)) \in \mathcal{C}(\mathbf{R}^d).$$

Hence the mapping $\Omega \ni \omega \mapsto \tilde{\mathcal{S}}(\omega) \in \mathcal{C}(\mathbf{R}^d)$ is $(\mathcal{F}, \mathcal{B}(\mathcal{C}(\mathbf{R}^d)))$-measurable.

A similar argument yields the measurability of $\Omega \ni \omega \mapsto \tilde{\mathcal{U}}(\omega) \in \mathcal{C}(\mathbf{R}^d)$. This completes the proof of assertion (h) of the theorem.

If $\overset{\circ}{F}$ (resp. $F$) satisfy Hypothesis $(\text{ST}(k,\delta))$ (resp. $(\text{IT}(k,\delta))$) for every $k \geq 1$ and $\delta \in (0,1]$, then a simple adaptation of the argument in [Ru.1], Section (5.3) (p. 297) gives a sure event in $\mathcal{F}$, also denoted by $\Omega^*$ such that $\tilde{\mathcal{S}}(\omega), \tilde{\mathcal{U}}(\omega)$ are $C^\infty$ for all $\omega \in \Omega^*$. This completes the proof of Theorem 3.1. $\square$

**Global Stable and Unstable Sets.**

We will conclude this section by a discussion of global stable and unstable sets for the sde's (S) and (I). Assume all the conditions of Theorem 3.1. Define the set

$$\tilde{\mathcal{S}}_g(\omega) := \{x \in \mathbf{R}^d : \limsup_{t \to \infty} \frac{1}{t} \log |\phi(t, x, \omega) - Y(\theta(t, \omega))| \leq \lambda_{i_0}\}$$

for each $\omega \in \Omega^*$. The family $\tilde{\mathcal{S}}_g(\omega), \omega \in \Omega^*$, is clearly invariant under $\phi$; that is

$$\phi(t, \cdot, \omega)(\tilde{\mathcal{S}}_g(\omega)) = \tilde{\mathcal{S}}_g(\theta(t, \omega))$$

for all $t \in \mathbf{R}$ and all $\omega \in \Omega^*$.

Using induction, we may define the family $\{\tilde{\mathcal{S}}^n(\omega)\}_{n=0}^\infty$ of $C^{k,\epsilon}$ stable submanifolds as follows:

$$\tilde{\mathcal{S}}^0(\omega) := \tilde{\mathcal{S}}(\omega)$$

$$\tilde{\mathcal{S}}^n(\omega) := \begin{cases} \phi(-n, \cdot, \omega)(\tilde{\mathcal{S}}(\theta(n, \omega))), & \text{if } \tilde{\mathcal{S}}^{n-1}(\omega) \subseteq \phi(-n, \cdot, \omega)(\tilde{\mathcal{S}}(\theta(n, \omega))) \\ \tilde{\mathcal{S}}^{n-1}(\omega), & \text{otherwise} \end{cases}$$

for $n \geq 1$. In the above definition, $\tilde{\mathcal{S}}(\omega)$ refers to the stable manifolds constructed in the proof of Theorem 3.1. Note that $\tilde{\mathcal{S}}^n(\omega) \subseteq \tilde{\mathcal{S}}^{n+1}(\omega)$ for all $n \geq 0$. Furthermore, the global stable set $\tilde{\mathcal{S}}_g(\omega)$ is given by

$$\tilde{\mathcal{S}}_g(\omega) = \bigcup_{n=1}^\infty \tilde{\mathcal{S}}^n(\omega), \quad \omega \in \Omega^*. \tag{38}$$

We will indicate a proof of (38). Fix any $\omega \in \Omega^*$. Then by asymptotic cocycle invariance of the stable manifolds, there is an a positive $l_0 := l_0(\omega)$ such that

$$\phi(l, \cdot, \omega)(\tilde{\mathcal{S}}(\omega)) \subseteq \tilde{\mathcal{S}}(\theta(l, \omega)) \tag{39}$$

for all integers $l \geq l_0$. The inclusion (39) follows from Remark (5.2)(c) in ([Ru.1], p. 297). In particular, and by the definition of $\tilde{\mathcal{S}}^n(\omega)$, it follows that $\tilde{\mathcal{S}}^n(\omega) = \phi(-n, \cdot, \omega)(\tilde{\mathcal{S}}(\theta(n, \omega)))$ for infinitely many integers $n > 0$. Now let $x \in \tilde{\mathcal{S}}_g(\omega)$. Then it is easy to see that $\phi(k, x, \omega) \in \tilde{\mathcal{S}}(\theta(k, \omega))$ for sufficiently large $k$. Fix such a $k$ and call it $k_0$. Then there exists $l \geq k_0$ such that $\phi(l, x, \omega) \in \tilde{\mathcal{S}}(\theta(l, \omega))$ and $\tilde{\mathcal{S}}^l(\omega) = \phi(-l, \cdot, \omega)(\tilde{\mathcal{S}}(\theta(l, \omega)))$. Hence $x \in \tilde{\mathcal{S}}^l(\omega)$ and therefore $x \in \bigcup_{n=1}^\infty \tilde{\mathcal{S}}^n(\omega)$. Conversely, let $x$ belong to the set on the right-hand-side of (38). Then by definition of the $\tilde{\mathcal{S}}^n(\omega)$, there exists $k$ such that $x \in \phi(-k, \cdot, \omega)(\tilde{\mathcal{S}}(\theta(k, \omega)))$. By Theorem 3.1(a), this implies that $\limsup_{t \to \infty} \frac{1}{t} \log |\phi(t, x, \omega) - Y(\theta(t, \omega))| \leq \lambda_{i_0}$. Hence $x \in \tilde{\mathcal{S}}_g(\omega)$, and the proof of (38) is complete.

Similar remarks hold for the global unstable set

$$\tilde{\mathcal{U}}_g(\omega) := \{x \in \mathbf{R}^d : \limsup_{t \to \infty} \frac{1}{t} \log |\phi(-t, x, \omega) - Y(\theta(-t, \omega))| \leq \lambda_{i_0-1}\}.$$

The above considerations also show that $\tilde{\mathcal{S}}_g$ and $\tilde{\mathcal{U}}_g$ are $C^{k,\epsilon}$ manifolds which are immersed (but not in general imbedded) in $\mathbf{R}^d$.

## 4. Appendix. The Substitution Rule.

In this appendix, we will establish some results that are aimed towards showing that certain extensions of the Itô integral and the Stratonovich integral are stable under random substitutions.

Throughout this appendix, $F : \mathbf{R} \times \mathbf{R}^l \times \Omega \to \mathbf{R}^d$ is a continuous spatial semimartingale based on a filtered probability space such that $F(0, x) = 0$ for all $x \in \mathbf{R}^l$. (Note that the helix property is not needed in this section.) We shall use the notation in Sections 1 and 2. Decompose $F$ as

$$F(t, x) = V(t, x) + M(t, x), \quad t \geq 0, \quad x \in \mathbf{R}^l \tag{1}$$

where $V(\cdot, x) := (V^1(\cdot, x), \cdots, V^d(\cdot, x))$ is a continuous bounded variation process and $M(\cdot, x) := (M^1(\cdot, x), \cdots, M^d(\cdot, x))$, $x \in \mathbf{R}^l$, is a continuous spatial local martingale such that $M(0, x) = V(0, x) = 0$ for all $x \in \mathbf{R}^l$ and all $\omega \in \Omega$.

We now introduce a definition of the Itô and the Stratonovich integral with respect to integrands that are possibly anticipating.

Let $f : [0, \infty) \times \Omega \to \mathbf{R}^l$ be a measurable process with continuous sample paths. Take any sequence of partitions $\pi_n := \{0 = t_0^n < t_1^n < \cdots < t_n^n\}$ of $[0, \infty)$. Suppose $\lim_{n \to \infty} t_n^n = \infty$ and $\lim_{n \to \infty} \max\{t_i^n - t_{i-1}^n, i = 1, \cdots, n\} = 0$. Define the sequence

$$I_n(t) := \sum_{k=0}^{n-1} [M(t_{k+1}^n \wedge t, f(t_k^n)) - M(t_k^n \wedge t, f(t_k^n))], \quad n \geq 1, t \geq 0. \tag{2}$$

If, in addition, $M$ is a $C^1$ spatial local martingale, then define

$$S_n(t) := I_n(t) + C_n(t), \quad n \geq 1, t \geq 0, \tag{3}$$

where

$$C_n(t) := \frac{1}{2} \sum_{k=0}^{n-1} [D_2 M(t_{k+1}^n \wedge t, f(t_k^n)) - D_2 M(t_k^n \wedge t, f(t_k^n))][f(t_{k+1}^n \wedge t) - f(t_k^n \wedge t)].$$

**Definition 4.1.**
(i) Define the Itô integral of $f$ with respect to the continuous spatial local martingale $M$ by

$$\int_0^T M(dt, f(t)) := \lim_{n \to \infty} I_n(T) \tag{4}$$

when the limit in probability exists uniformly on compact subsets of $[0, \infty)$ for any sequence of partitions as above.

(ii) Define the Stratonovich integral of $f$ with respect to the $C^1$ local martingale $M$ by

$$\int_0^T M(\circ dt, f(t)) := \lim_{n \to \infty} S_n(T) \tag{5}$$

when the limit exists uniformly on compact subsets of $[0, \infty)$ in probability for any sequence of partitions as above.

(iii) If $F$ is a spatial semimartingale given by (1), define the Stratonovich integral of $f$ with respect to $F$ by

$$\int_0^T F(\circ dt, f(t)) := \int_0^T V(dt, f(t)) + \int_0^T M(\circ dt, f(t)). \tag{6}$$

provided the right side of (6) is defined. The Itô integral is defined analogously (without the circle).

Note that our definitions of the Itô and the Stratonovich integral agree with the classical ones when the integrand process $f$ is a continuous semimartingale. For the Stratonovich integral this follows from ([Ku], Theorem 3.2.5, p. 86). As will be clear from the sequel, the computations become simpler under Definition 4.1 than if we had directly extended Kunita's definition to the non-adapted case. We remark that our definition of the Itô integral does not always coincide with the well-known Skorohod integral even if both are defined.

In the following theorem, $B_c^{k,\delta}$ denotes the class of all $C^k$ spatial semimartingales such that for any $T > 0$, any $p \geq 1$ and any compact subset $K$ of $\mathbf{R}^l$ (or $\mathbf{R}^m$) the $p$-th moment of the $(k+\delta)$-norms of the characteristics restricted to $K$ are uniformly bounded on $[0, T]$. Observe that $B_{ub}^{k,\delta} \subset B_c^{k,\delta}$. The flow generated by an Itô equation driven by $F \in B_{ub}^{0,1}$ is always in $B_c^{0,1}$ but generally not in $B_{ub}^{0,1}$ (see Theorem 4.2).

We now state the substitution rule.

**Theorem 4.1.**
*Fix $\delta \in (0,1)$ and let $F(t,y) = M(t,y) + V(t,y) \in \mathbf{R}^d$, $y \in \mathbf{R}^l$, be a spatial semimartingale of class $B_c^{0,\delta}$ such that $M(0,y) = V(0,y) = 0$ for all $y \in \mathbf{R}^l$. Further let $f : [0, \infty) \times \mathbf{R}^m \times \Omega \to \mathbf{R}^l$ be a continuous spatial semimartingale such that for any compact subset $K$ of $\mathbf{R}^m$, any $T > 0$ and any $p > 1$, there exists a constant $c$ such that $E(|f(t,x) - f(t,y)|^p) \leq c|x-y|^{\delta p}$ for all $x, y \in K$ and all $0 \leq t \leq T$.*

*Then there is a modification of the Itô integral such that for any $\mathcal{F}$-measurable random variable $Y : \Omega \to \mathbf{R}^m$, one has a.s.*

$$\int_0^T F(dt, f(t,x))\Big|_{x=Y} = \int_0^T F(dt, f(t,Y)) \tag{7}$$

*for all $T > 0$.*

*If, moreover, $M$ is of class $B_c^{1,1}$ and $f \in B_c^{0,\delta}$, then we also have*

$$\int_0^T F(\circ dt, f(t,x))\Big|_{x=Y} = \int_0^T F(\circ dt, f(t,Y)) \quad a.s. \tag{8}$$

*for all $T > 0$.*

For Brownian linear integrators, a similar result is given in ([N-P], Propositions 7.7, 7.8), ([A-I], Theorem 2, Corollary 1) and ([Nu] (Theorem 5.3.3). In order to prove Theorem 4.1, we will adopt the approach in ([Nu]). The essence of the argument is to replace $f$ by $f(t, x)$ on the right-hand-side of (2), substitute $x = Y(\omega)$ in each finite sum in (2), and then pass to the limit in probability in order to get (7) and (8).

Note that the substitution rule holds trivially in the bounded variation integral on the right-hand-side of (6). Hence in all subsequent computations, we can and will assume that $V \equiv 0$ and $F = M$.

The proof of Theorem 4.1 turns on the following lemma.

**Lemma 4.1.**

*Let $\{S_n(x), x \in \mathbf{R}^m\}$, $n \geq 1$ be a sequence of (jointly) measurable random fields taking values in a complete separable metric space $(E, \rho)$ such that*

$$\lim_{n \to \infty} S_n(x) = S(x)$$

*in probability, where $\{S(x), x \in \mathbf{R}^m\}$ is a random field. Assume that there exist positive constants $p \geq 1$, $\alpha > m$, $C = C(T, K, p)$ such that whenever $K > 0$, and $|x|, |x'| \leq K$, one has*

$$E[\rho(S_n(x), S_n(x'))^p] \leq C|x - x'|^\alpha$$

*for all $n \geq 1$. Then the random fields $S$, $S_n$, $n \geq 1$, have continuous modifications (denoted by the same symbols). For any such modifications and any random variable $Y : \Omega \to \mathbf{R}^m$, one has*

$$\lim_{n \to \infty} S_n(Y) = S(Y)$$

*in probability.*

The proof is given in ([Nu], Lemma 5.3.1) in the special case where $E = \mathbf{R}^d$ but the argument therein carries over to our case without change. Observe that the conditions of the Lemma imply that $S_n(\cdot)$ converges to $S(\cdot)$ uniformly on compact subsets of $\mathbf{R}^m$ in probability (which is the reason why the substitution property holds).

*Proof of Theorem 4.1.*

Assume, without loss of generality, that $F = M$, a local martingale. Let us first assume in addition that $M$ and $f$ are bounded uniformly in $(t, x, \omega)$ on compact subsets of $[0, \infty) \times \mathbf{R}^l$ (resp. $[0, \infty) \times \mathbf{R}^m$). For a given sequence of partitions $\pi_n := \{0 = t_0^n < t_1^n < \cdots < t_n^n\}$ of $[0, \infty)$ as in Definition 4.1 define the sequence of random fields $I_n(t, x)$, $x \in \mathbf{R}^m, n \geq 1$, by

$$I_n(t, x) := \sum_{k=0}^{n-1} [M(t_{k+1}^n \wedge t, f(t_k^n, x)) - M(t_k^n \wedge t, f(t_k^n, x))], \quad n \geq 1, t \geq 0. \quad (9)$$

We want to check the assumptions of Lemma 4.1 for the sequence $I_n$ taking values in the space of $\mathbf{R}^d$-valued continuous functions on $[0, \infty)$. Write $I_n(t, x) := (I_n^1(t, x), \cdots, I_n^d(t, x))$. It is enough to show that for every compact subset $K$ of $\mathbf{R}^m$, every $T > 0$ and every $j \in \{1, 2, \cdots, d\}$, we have

$$E[\sup_{0 \leq s \leq T} |I_n^j(s, x) - I_n^j(s, y)|^p] \leq c|x - y|^\alpha$$

for some $p > 1$, some $\alpha > m$, some $c > 0$ and all $x, y \in K$.

Fix a compact subset $K$ of $\mathbf{R}^m$, $T > 0$, $j \in \{1, \cdots, d\}$, $p > 1$, and abbreviate $u_k = t_k^n \wedge T$. Using the Burkholder-Davis-Gundy inequality, we find constants $c_1$, $c_2$ and $c_3$ (independent of $n$) such that for all $x, y \in K$, one has

$$E[\sup_{0 \leq s \leq T} |I_n^j(s, x) - I_n^j(s, y)|^p]$$

$$\leq c_1 E \left[ \sum_{k=0}^{n-1} \left\{ M^j(u_{k+1}, f(u_k, x)) - M^j(u_{k+1}, f(u_k, y)) - \right. \right.$$
$$\left. \left. + -M^j(u_k, f(u_k, x)) + M^j(u_k, f(u_k, y)) \right\}^2 \right]^{p/2}$$

$$\leq c_1 \left[ \sum_{k=0}^{n-1} \left\{ E|M^j(u_{k+1}, f(u_k, x)) - M^j(u_{k+1}, f(u_k, y)) \right. \right.$$
$$\left. \left. + -M^j(u_k, f(u_k, x)) + M^j(u_k, f(u_k, y))|^p \right\}^{2/p} \right]^{p/2}$$

$$\leq c_2 \left[ \sum_{k=0}^{n-1} \left\{ E(\int_{u_k}^{u_{k+1}} a_{jj}(s, f(u_k, x), f(u_k, x)) - 2a_{jj}(s, f(u_k, x), f(u_k, y)) \right. \right.$$
$$\left. \left. + a_{jj}(s, f(u_k, y), f(u_k, y)) ds)^{p/2} \right\}^{2/p} \right]^{p/2}$$

$$\leq c_3 \left[ \sum_{k=0}^{n-1} (u_{k+1} - u_k)(E(|f(u_k, x) - f(u_k, y)|^{2\delta p}))^{1/p} \right]^{p/2},$$

where we have used Hölder's inequality, the boundedness of $f$ and the fact that $M \in B_c^{0,\delta}$ to obtain the last inequality.

Inserting the moment estimate on $f$ in Theorem 4.1, we see that for each $p \geq 1, T > 0$ and each compact subset $K$ of $\mathbf{R}^m$, there exists a constant $c_4$ such that

$$E[\sup_{0 \leq s \leq T} |I_n(s, x) - I_n(s, y)|^p] \leq c_4 |x - y|^{\delta^2 p}.$$

for all $n$ and all $x, y \in K$. Now take $p$ sufficiently large so that $\delta^2 p > m$. Therefore the substitution formula follows from Lemma 4.1 in the Itô case under the additional constraint that $M$ and $f$ are bounded. For general $M$ and $f$, we get the uniform convergence $I_n(t,x) \to I(t,x)$ in probability on compacts of $[0,\infty) \times \mathbf{R}^m$, by a straightforward localization argument.

To show (8) we assume first that $M$, $D_2 M$ and $f$ are uniformly bounded in $(t,x,\omega)$ on compact subsets. Let

$$C_n(t,x) = \frac{1}{2} \sum_{k=0}^{n-1} [D_2 M(t_{k+1}^n \wedge t, f(t_k^n, x)) - D_2 M(t_k^n \wedge t, f(t_k^n, x))][f(t_{k+1}^n \wedge t, x) - f(t_k^n \wedge t, x)].$$

To apply Lemma 4.1, we will show that for every compact subset $K$ of $\mathbf{R}^m$ and every $T > 0$ there exist $p > 1$, $\alpha > m$ and $c_5 > 0$ such that

$$E \sup_{0 \le s \le T} |C_n(s,x) - C_n(s,y)|^p \le c_5 |x - y|^\alpha$$

for all $n \in \mathbf{N}$ and all $x, y \in K$.

We will use the following abbreviations (suppressing the dependence on $n$):

$$A_k(t,x) = D_2 M(t_{k+1}^n \wedge t, f(t_k^n, x)) - D_2 M(t_k^n \wedge t, f(t_k^n, x))$$

and

$$B_k(t,x) = f(t_{k+1}^n \wedge t, x) - f(t_k^n \wedge t, x).$$

Then we get for $p \ge 1$ and all $x, y \in K$

$$E \sup_{0 \le t \le T} |2C_n(t,x) - 2C_n(t,y)|^p \qquad (10)$$

$$\le 2^p (E \sup_{0 \le t \le T} |\sum_{k=0}^{n-1} (A_k(t,x) - A_k(t,y)) B_k(t,x)|^p + E \sup_{0 \le t \le T} |\sum_{k=0}^{n-1} A_k(t,y)(B_k(t,x) - B_k(t,y))|^p).$$

Using Hölder's inequality, we get

$$E \sup_{0 \le t \le T} |\sum_{k=0}^{n-1} (A_k(t,x) - A_k(t,y)) B_k(t,x)|^p$$

$$\le (\sum_{k=0}^{n-1} (E((\sup_{0 \le t \le T} |A_k(t,x) - A_k(t,y)|^p)(\sup_{0 \le t \le T} |B_k(t,x)|^p)))^{1/p})^p$$

$$\le (\sum_{k=0}^{n-1} (E \sup_{0 \le t \le T} |A_k(t,x) - A_k(t,y)|^{2p})^{1/2p} (E \sup_{0 \le t \le T} |B_k(t,x)|^{2p})^{1/2p})^p. \qquad (11)$$

Using the inequality of Burkholder, Davis and Gundy and [Ku], Theorem 3.1.2 which allows us to interchange spatial derivatives and the quadratic variation there exist constants $c_6$, $c_7$ and $c_8$ (independent of $k$ and $n$) such that

$$(E \sup_{0 \le t \le T} |A_k(t,x) - A_k(t,y)|^{2p})^{1/2p} \le c_6 (u_{k+1} - u_k)^{1/2} |x - y|^{\delta^2} \qquad (12)$$

(this is derived just like the first part of the proof) and

$$(E \sup_{0 \le t \le T} (|B_k(t,x)|^{2p}))^{1/2p} \le c_7(u_{k+1} - u_k)^{1/2} + c_8(u_{k+1} - u_k). \tag{13}$$

Since $(u_{k+1} - u_k) \le T^{1/2}(u_{k+1} - u_k)^{1/2}$ we can in fact delete the term $c_8(u_{k+1} - u_k)$ in (13) by increasing $c_7$ accordingly.

Similarly, for $p \ge 1$ there exists some constant $c_9$ such that for all $x, y \in K$, we have

$$E \sup_{0 \le t \le T} |\sum_{k=0}^{n-1} A_k(t,y)(B_k(t,x) - B_k(t,y))|^p$$

$$\le \left[\sum_{k=0}^{n-1} (E \sup_{0 \le t \le T} |A_k(t,y)|^{2p})^{1/2p} (E \sup_{0 \le t \le T} |B_k(t,x) - B_k(t,y)|^{2p})^{1/2p}\right]^p$$

$$\le c_9 (\sum_{k=0}^{n-1} ((u_{k+1} - u_k)^{1/2} |x-y|^\delta (u_{k+1} - u_k)^{1/2}))^p \le c_9 T^p |x-y|^{\delta p}. \tag{14}$$

Inserting (12) and (13) into (11) and then (11) and (14) into (10) we see that for all $p \ge 1$ there exists a constant $c_5$ such that for all $x, y \in K$ we have

$$E \sup_{0 \le t \le T} |2C_n(t,x) - 2C_n(t,y)|^p \le c_5 |x-y|^{\delta p}$$

as desired. The case of general $M$ and $f$ is again easily obtained by localization. This proves Theorem 4.1. □

The next result allows us to start the solution of a sde and its linearization at any random (possibly anticipating) initial state. In all parts of the following theorem the stochastic flow associated with an Itô equation driven by $F$ or a Stratonovich equation driven by $\overset{\circ}{F}$ will be denoted by $\phi$. Recall that $\phi(t,x) := \phi_{0t}(x)$.

**Theorem 4.2.**

Suppose $Y: \Omega \to \mathbf{R}^d$ is any $\mathcal{F}$-measurable random variable and let $\delta \in (0,1]$.

(i) Let $F$ be a spatial (forward) semimartingale in $B_{ub}^{0,1}$. Then $\phi(t,Y), t \ge 0$, is a solution of the anticipating Itô sde

$$\left. \begin{aligned} d\phi(t,Y) &= F(dt, \phi(t,Y)), \quad t > 0 \\ \phi(0,Y) &= Y. \end{aligned} \right\} \tag{I'}$$

If the spatial (forward) semimartingale $\overset{\circ}{F}$ is in $(B_{ub}^{2,\delta}, B_{ub}^{1,\delta})$, then $\phi(t,Y), t \ge 0$, is a solution of the anticipating Stratonovich sde

$$\left. \begin{aligned} d\phi(t,Y) &= \overset{\circ}{F}(\circ dt, \phi(t,Y)), \quad t > 0 \\ \phi(0,Y) &= Y. \end{aligned} \right\} \tag{SII}$$

(ii) Assume that $\overset{\circ}{F}$ is a (forward) spatial semimartingale of class $(B_{ub}^{3,\delta}, B_{ub}^{2,\delta})$. Then the (possibly non-adapted) process $y(t, \omega) := D_2\phi(t, Y(\omega), \omega)$, $t \geq 0$, satisfies the Stratonovich linearized sde

$$\left.\begin{aligned} dy(t) &= D_2\overset{\circ}{F}(\circ dt, \phi(t, Y))y(t), \quad t > 0, \\ y(0) &= I \in L(\mathbf{R}^d). \end{aligned}\right\} \quad (SIII)$$

A similar result is true in the Itô case.

(iii) Let $\overset{\circ}{F}$ be a spatial backward semimartingale of class $(B_{ub}^{2,\delta}, B_{ub}^{1,\delta})$. Then $\phi(t, Y), t \leq 0$, is a solution of the backward Stratonovich sde

$$\left.\begin{aligned} d\phi(t, Y) &= -\overset{\circ}{F}(\circ \hat{d}t, \phi(t, Y)), \quad t < 0 \\ \phi(0, Y) &= Y. \end{aligned}\right\} \quad (SII^-)$$

(iv) Assume that $\overset{\circ}{F}$ is a spatial backward semimartingale of class $(B_{ub}^{3,\delta}, B_{ub}^{2,\delta})$. Then the process $y(t, \omega) := D_2\phi(t, Y(\omega), \omega)$, $t \leq 0$, satisfies the backward Stratonovich linearized sde

$$\left.\begin{aligned} d\hat{y}(t) &= -D_2\overset{\circ}{F}(\circ \hat{d}t, \phi(t, Y))\hat{y}(t), \quad t < 0, \\ \hat{y}(0) &= I \in L(\mathbf{R}^d). \end{aligned}\right\} \quad (SIII^-)$$

*Proof.*

Let $F$ be in $B_{ub}^{0,1}$ and define $f(t, x) := \phi(t, x)$. Then the moment estimate for $f$ in Theorem 4.1 is satisfied (with $\delta = 1$), thanks to [Ku], Lemma 4.5.6. Therefore (I$'$) follows.

Next suppose that $\overset{\circ}{F}$ is in $(B_{ub}^{2,\delta}, B_{ub}^{1,\delta})$, so in particular the local martingale part is in $B_{ub}^{1,1}$. By [Ku], Theorem 3.4.7 (or our Proposition 2.1) we know that $\phi$ is also generated by an Itô equation which is driven by a semimartingale $F$ with local characteristics of class $(B_{ub}^{2,\delta}, B_{ub}^{1,\delta})$. Observe that $f := \phi_{0\cdot} \in B_c^{0,1}$ because, for every compact subset $K \subset \mathbf{R}^d$, every $T > 0$ and every $p \geq 1$, we have

$$\sup_{0 \leq s \leq T} E\left[\sup_{x \in K} |\phi_{0s}(x)|^p\right] < \infty, \tag{15}$$

and

$$\sup_{0 \leq s \leq T} E\left[\sup_{x,y \in K, x \neq y} \left(\frac{|\phi_{0s}(x) - \phi_{0s}(y)|}{|x - y|}\right)^p\right] < \infty. \tag{16}$$

The estimates (15) and (16) follow from Theorem 2.1(v). This proves part (i).

We next proceed to prove assertion (ii) of the theorem. To do this we will reduce the problem to a system of sde's which satisfies the hypotheses of part (i) of the theorem. Define the spatial semimartingales

$$z(t, x, v) := (\phi(t, x), D_2\phi(t, x)(v))$$

$$G(t, x, v) := (\overset{\circ}{F}(t, x), D_2\overset{\circ}{F}(t, x)(v))$$

for all $(x, v) \in \mathbf{R}^d \times \mathbf{R}^d$, $t > 0$. Then the sde's (S) and its linearization

$$d[D_2\phi(t, x)(v)] = D_2\overset{\circ}{F}(\circ dt, \phi(t, x))D_2\phi(t, x)(v), \quad t > 0,$$

$$D_2\phi(0, x)(v) = v \in \mathbf{R}^d$$

viewed as a coupled pair, are equivalent to the sde

$$\left.\begin{array}{l} dz(t, x, v) = G(\circ dt, z(t, x, v)), \quad t > 0 \\ z(0, x, v) = (x, v) \in \mathbf{R}^d \times \mathbf{R}^d. \end{array}\right\} \quad (SIV)$$

By hypothesis, $\overset{\circ}{F}$ has local characteristics of class $(B_{ub}^{3,\delta}, B_{ub}^{2,\delta})$. We claim that $G$ has local characteristics of class $(B_{ub}^{2,\delta}, B_{ub}^{1,\delta})$ (in $\mathbf{R}^{2d}$). To see this, use coordinates

$$x := (x_1, \cdots, x_d), \quad x' := (x'_1, \cdots, x'_d), \quad v := (v_1, \cdots, v_d), \quad v' := (v'_1, \cdots, v'_d),$$

and observe that the semimartingale $\{D_2\overset{\circ}{F}(t, x)(v) : (x, v) \in \mathbf{R}^d \times \mathbf{R}^d, t \geq 0\}$ has local characteristics $\{\tilde{a}^{k,l}(t, (x, v), (x', v')) : (x, v), (x', v') \in \mathbf{R}^d \times \mathbf{R}^d, t \geq 0, 1 \leq k, l \leq d\}$, $\{\tilde{b}^k(t, (x, v)) : (x, v) \in \mathbf{R}^d \times \mathbf{R}^d, t \geq 0, 1 \leq k \leq d\}$ given by

$$\tilde{a}^{k,l}(t, (x, v), (x', v')) = \sum_{i,j=1}^{d} \frac{\partial^2}{\partial x_i \partial x'_j} a^{k,l}(t, x, x') v_i v'_j$$

$$\tilde{b}^k(t, (x, v)) = \sum_{i=1}^{d} \frac{\partial}{\partial x_i} b^k(t, x) v_i.$$

From these relations, our claim follows. By the first part of the proof, we can substitute $x = Y(\omega)$ in (SIV) and keep $v \in \mathbf{R}^d$ arbitrary but fixed (non-random). This gives (SIII) (and (SII)). Hence assertion (ii) of the theorem holds.

The proofs of assertions (iii) and (iv) are similar to those of (i) and (ii). This completes the proof of the theorem. □

In the case of Brownian linear integrators, a version of Theorem 4.2 (i) is given in ([M-N-S], Theorem 3.1, p. 1920) under somewhat more restrictive hypotheses. In this case too, similar results to Theorem 4.2 (i), (ii) appear in ([A-I], Theorems 4,5) and ([Nu], Theorems 5.3.4, 6.1.1).

**Acknowledgement.**
The authors are grateful for comments made by Ludwig Arnold on an earlier version of the manuscript.

# References

[A]      Arnold, L., *Random Dynamical Systems*, Springer-Verlag (To appear).

[A-S]    Arnold, L., and Scheutzow, M. K. R., Perfect cocycles through stochastic differential equations, *Probab. Th. Rel. Fields,* 101, (1995), 65-88.

[A-I]    Arnold, L., and Imkeller, P., Stratonovich calculus with spatial parameters and anticipative problems in multiplicative ergodic theory, *Stoch. Proc. Appl.* 62 (1996), 19-54.

[Bo]     Boxler, P., A stochastic version of center manifold theory, *Probab. Th. Rel. Fields*, 83 (1989), 509-545.

[C]      Carverhill, A., Flows of stochastic dynamical systems: Ergodic theory, *Stochastics*, 14 (1985), 273-317.

[I-S]    Imkeller, P., and Scheutzow, M. K. R., On the spatial asymptotic behaviour of stochastic flows in Euclidean space (preprint)(1997).

[I-W]    Ikeda, N., and Watanabe, S., *Stochastic Differential Equations and Diffusion Processes,* Second Edition, North-Holland-Kodansha (1989).

[Ku]     Kunita, H., *Stochastic Flows and Stochastic Differential Equations*, Cambridge University Press, Cambridge, New York, Melbourne, Sydney (1990).

[Mo.1]   Mohammed, S.-E. A., The Lyapunov spectrum and stable manifolds for stochastic linear delay equations, *Stochastics and Stochastic Reports*, Vol. 29 (1990), 89-131.

[Mo.2]   Mohammed, S.-E. A., Lyapunov exponents and stochastic flows of linear and affine hereditary Systems, *Diffusion Processes and Related Problems in Analysis*, Vol. II, edited by Mark Pinsky and Volker Wihstutz, Birkhauser (1992), 141-169.

[M-N-S]  Millet, A., Nualart, D., and Sanz, M., Large deviations for a class of anticipating stochastic differential equations, *The Annals of Probability*, 20 (1992), 1902-1931.

[M-S.1]  Mohammed, S.-E. A., and Scheutzow, M. K. R., Lyapunov exponents of linear stochastic functional differential equations driven by semimartingales, Part I: The multiplicative ergodic theory, *Ann. Inst. Henri Poincaré, Probabilités et Statistiques,* Vol. 32, 1, (1996), 69-105. pp. 43.

[M-S.2]  Mohammed, S.-E. A., and Scheutzow, M. K. R., Spatial estimates for stochastic flows in Euclidean space, to appear in *The Annals of Probability.*

[Nu]     Nualart, D., Analysis on Wiener space and anticipating stochastic calculus (to appear in) *St. Flour Notes.*

STABLE MANIFOLD THEOREM FOR SDE'S 37[N-P]   Nualart, D., and Pardoux, E., Stochastic calculus with anticipating integrands, Analysis on Wiener space and anticipating stochastic calculus, *Probab. Th. Rel. Fields*, 78 (1988), 535-581.

[O]     Oseledec, V. I., A multiplicative ergodic theorem. Lyapunov characteristic numbers for dynamical systems, *Trudy Moskov. Mat. Obšč.* 19 (1968), 179-210. English transl. *Trans. Moscow Math. Soc.* 19 (1968), 197-221.

[O-P]   Ocone, D., and Pardoux, E., A generalized Itô-Ventzell formula. Application to a class of anticipating stochastic differential equations, *Ann. Inst. Henri Poincaré, Probabilités et Statistiques,* Vol. 25, no. 1 (1989), 39-71.

[Pr.1]  Protter, Ph. E., Semimartingales and measure-preserving flows, *Ann. Inst. Henri Poincaré, Probabilités et Statistiques,* vol. 22, (1986), 127-147.

[Pr.2]  Protter, Ph. E., *Stochastic Integration and Stochastic Differential Equations: A New Approach*, Springer (1990).

[Ru.1]  Ruelle, D., Ergodic theory of differentiable dynamical systems, *Publ. Math. Inst. Hautes Etud. Sci.* (1979), 275-306.

[Ru.2]  Ruelle, D., Characteristic exponents and invariant manifolds in Hilbert space, *Annals of Mathematics 115* (1982), 243–290.

[Sc]    Scheutzow, M. K. R., On the perfection of crude cocycles, *Random and Computational Dynamics, 4,* (1996), 235-255.

[Wa]    Wanner, T., Linearization of random dynamical systems, in *Dynamics Reported,* vol. 4, edited by U. Kirchgraber and H.O. Walther, Springer (1995), 203-269.